\documentclass{article}

\PassOptionsToPackage{%
	style	= alphabetic%
}{biblatex}

\usepackage{geometry}
\usepackage[hidelinks]{hyperref}
\hypersetup{
	colorlinks,
	linkcolor={blue},
	citecolor={ForestGreen},
	urlcolor={blue}
}

\usepackage{graphicx}
\usepackage[font = normal, labelfont = sf]{subcaption}
\graphicspath{{TCSI_Figures/}}
\usepackage[font = normal, labelfont = sf, labelsep = period, singlelinecheck = false]{caption}

\usepackage[UKenglish]{babel}
\usepackage{amsmath, amsthm, amssymb, mathrsfs, bm, mathtools}
\usepackage{wasysym}
\allowdisplaybreaks

\usepackage{titlesec}
\usepackage{xifthen}
\usepackage{xspace}

\usepackage{xparse}

\usepackage{titlesec}
\titleformat{\title}
{\sffamily \Huge \bfseries \scshape \boldmath}{}{}{}
\titleformat{\section}
{\sffamily \Large \bfseries \scshape \boldmath}{\thesection}{1em}{}
\titleformat{\subsection}
{\sffamily \large \bfseries \scshape \boldmath}{\thesubsection}{1em}{}
\titleformat{\subsubsection}
{\sffamily \normalsize \bfseries \scshape \boldmath}{\thesubsubsection}{1em}{}
\titleformat{\paragraph}[runin]
{\sffamily \normalsize \bfseries \scshape \boldmath}{\theparagraph}{1em}{}
\titleformat{\subparagraph}[runin]
{\sffamily \normalsize \bfseries \scshape \boldmath}{\thesubparagraph}{1em}{}

\usepackage{enumitem}

\setlist[enumerate]{%
	topsep	= \medskipamount,
	itemsep = 0pt,
	label	= \textup{(\roman*)}
}
\setlist[itemize]{%
	topsep	= \medskipamount,
	itemsep	= 0pt,
	label	= \bcdot
}
\setlist[description]{%
	topsep	= \smallskipamount,		
	itemsep	= \smallskipamount,		
	font	= {\mdseries\itshape},	
}

\newlist{ranlist}{enumerate}{3}
\setlist[ranlist,1]{
	label=(\roman*) }
\setlist[ranlist,2]{
	label=(\textit{\alph*}),
	ref=(\roman{ranlisti}.\textit{\alph*})	}
\setlist[ranlist,3]{
	label=\arabic*.,
	ref=(\roman{ranlisti}.\textit{\alph{ranlistii}}.\arabic*)	}

\newcommand{\negphantom}[1]{\ifmmode\settowidth{\dimen0}{$#1$}\else\settowidth{\dimen0}{#1}\fi\hspace*{-\dimen0}}

\newenvironment{center-small}
	{\par\centering\smallskip}
	{\par\smallskip\noindent}

\usepackage[dvipsnames,hyperref]{xcolor}

\newcommand{\Quad}[1]{
	\mathchoice
	{\quad\text{#1}\quad}
	{\text{ #1 }}
	{\text{ #1 }}
	{\text{ #1 }}
}

\newcommand{\Qand}{\Quad{and}}
\newcommand{\Qas}{\Quad{as}}

\newcommand{\Qfor}{\Quad{for}}
\newcommand{\Qforall}{\Quad{for all}}

\newcommand{\Qif}{\Quad{if}}

\newcommand{\Qwhere}{\Quad{where}}
\newcommand{\Qwhen}{\Quad{when}}

\newcommand{\one}  [1]{\bm1\{ #1 \}}

\let\originalexp\exp
\let\exp\relax
\DeclareRobustCommand{\exp} [1]{\originalexp(#1)}
\newcommand{\expb} [1]{\originalexp\bigl( #1 \bigr)}

\newcommand{\bcdot}{\ensuremath{\bm{\cdot}}}
\newcommand{\cq}{\coloneqq}

\newcommand{\abs}  [1]{| #1 |}
\newcommand{\absb} [1]{\bigl| #1 \bigr|}

\newcommand{\absBB}[1]{\Biggl| #1 \Biggr|}

\newcommand{\rbr} [1]{ ( #1 ) }
\newcommand{\rbb} [1]{\bigl( #1 \bigr)}

\newcommand{\rbbb}[1]{\biggl( #1 \biggr)}

\newcommand{\bra} [1]{ \{ #1 \} }
\newcommand{\brb} [1]{\bigl\{ #1 \bigr\}}

\newcommand{\midb}{\bigm|}

\newcommand{\Oh}  [1]{\mathcal{O}( #1 )}

\newcommand{\oh}  [1]{o( #1 )}

\newcommand{\pr}[2][]{
	\mathchoice
	{\ifthenelse{\isempty{#1}}
		{\mathbb{P}\bigl(#2\bigr)}
		{\mathbb{P}_{#1}\bigl(#2\bigr)}}
	{\ifthenelse{\isempty{#1}}
		{\mathbb{P}(#2)}
		{\mathbb{P}_{#1}(#2)}}
	{\ifthenelse{\isempty{#1}}
		{\mathbb{P}(#2)}
		{\mathbb{P}_{#1}(#2)}}
	{\ifthenelse{\isempty{#1}}
		{\mathbb{P}(#2)}
		{\mathbb{P}_{#1}(#2)}}
}
\newcommand{\prt}[2][]{
	\ifthenelse{\equal{}{#1}}
	{\mathbb{P}( #2 )}
	{\mathbb{P}_{#1}( #2 )}
}
\newcommand{\prb}[2][]{
	\ifthenelse{\equal{}{#1}}
	{\mathbb{P}\bigl( #2 \bigr)}
	{\mathbb{P}_{#1}\bigl( #2 \bigr)}
}
\newcommand{\prB}[2][]{
	\ifthenelse{\equal{}{#1}}
	{\mathbb{P}\Bigl( #2 \Bigr)}
	{\mathbb{P}_{#1}\Bigl( #2 \Bigr)}
}
\newcommand{\prbb}[2][]{
	\ifthenelse{\equal{}{#1}}
	{\mathbb{P}\biggl( #2 \biggr)}
	{\mathbb{P}_{#1}\biggl( #2 \biggr)}
}
\newcommand{\prBB}[2][]{
	\ifthenelse{\equal{}{#1}}
	{\mathbb{P}\Biggl( #2 \Biggr)}
	{\mathbb{P}_{#1}\Biggl( #2 \Biggr)}
}
\newcommand{\prs}[2][]{
	\ifthenelse{\equal{}{#1}}
	{\mathbb{P}\left( #2 \right)}
	{\mathbb{P}_{#1}\left( #2 \right)}
}

\newcommand{\ex}[2][]{
	\mathchoice
	{\ifthenelse{\isempty{#1}}
		{\mathbb{E}\bigl(#2\bigr)}
		{\mathbb{E}_{#1}\bigl(#2\bigr)}}
	{\ifthenelse{\isempty{#1}}
		{\mathbb{E}(#2)}
		{\mathbb{E}_{#1}(#2)}}
	{\ifthenelse{\isempty{#1}}
		{\mathbb{E}(#2)}
		{\mathbb{E}_{#1}(#2)}}
	{\ifthenelse{\isempty{#1}}
		{\mathbb{E}(#2)}
		{\mathbb{E}_{#1}(#2)}}
}
\newcommand{\ext}[2][]{
	\ifthenelse{\equal{}{#1}}
	{\mathbb{E}( #2 )}
	{\mathbb{E}_{#1}( #2 )}
}
\newcommand{\exb}[2][]{
	\ifthenelse{\equal{}{#1}}
	{\mathbb{E}\bigl( #2 \bigr)}
	{\mathbb{E}_{#1}\bigr( #2 \bigr)}
}
\newcommand{\exB}[2][]{
	\ifthenelse{\equal{}{#1}}
	{\mathbb{E}\Bigl( #2 \Bigr)}
	{\mathbb{E}_{#1}\Bigl( #2 \Bigr)}
}
\newcommand{\exbb}[2][]{
	\ifthenelse{\equal{}{#1}}
	{\mathbb{E}\biggl( #2 \biggr)}
	{\mathbb{E}_{#1}\biggl( #2 \biggr)}
}
\newcommand{\exBB}[2][]{
	\ifthenelse{\equal{}{#1}}
	{\mathbb{E}\Biggl( #2 \Biggr)}
	{\mathbb{E}_{#1}\Biggl( #2 \Biggr)}
}

\newcommand{\Varb}[2][]{
	\ifthenelse{\equal{}{#1}}
	{\mathbb{V}\textnormal{ar} \bigl(#2\bigr)}
	{\mathbb{V}\textnormal{ar}_{#1} \bigl(#2\bigr)}
}
\newcommand{\VAR}[2][]{
	\ifthenelse{\equal{}{#1}}
	{\textnormal{Var}(#2)}
	{\textnormal{Var}_{#1}(#2)}
}

\usepackage{calc}
\newlength{\halfplusheight}
\newcommand{\LIM}[1]{%
	\mathop{\raisebox{\halfplusheight}{\(\displaystyle\lim_{#1}\)}}%
}

\newcommand{\maxt}[1]{
	\mathchoice
	{\textstyle \max_{#1} \displaystyle}
	{\max_{#1}}
	{\max_{#1}}
	{\max_{#1}}
}

\newcommand{\LIMSUP}[1]{%
	\mathop{\raisebox{\halfplusheight}{\(\displaystyle\limsup_{#1}\)}}%
}

\newcommand{\binomt}[2]{
	\mathchoice
	{\textstyle \binom{#1}{#2} \displaystyle}
	{\binom{#1}{#2}}
	{\binom{#1}{#2}}
	{\binom{#1}{#2}}
}

\NewDocumentCommand{\sumt}{smo}{%
	\mathchoice%
	{\textstyle\sum_{#2}\IfBooleanT{#1}{^\star}\IfValueT{#3}{^{#3}}\displaystyle}%
	{\sum_{#2}\IfBooleanT{#1}{^\star}\IfValueT{#3}{^{#3}}}%
	{\sum_{#2}\IfBooleanT{#1}{^\star}\IfValueT{#3}{^{#3}}}%
	{\sum_{#2}\IfBooleanT{#1}{^\star}\IfValueT{#3}{^{#3}}}%
}

\NewDocumentCommand{\sumd}{smo}{%
	\sum_{#2}\IfBooleanTF{#1}{^\star}{\IfValueT{#3}{^{#3}}}%
}

\newcommand{\intt}[2][]{
	\mathchoice
	{\ifthenelse{\isempty{#1}}
		{\textstyle \int_{#2}      \displaystyle}
		{\textstyle \int_{#2}^{#1} \displaystyle}}
	{\ifthenelse{\isempty{#1}}
		{\int_{#2}}
		{\int_{#2}^{#1}}}
	{\ifthenelse{\isempty{#1}}
		{\int_{#2}}
		{\int_{#2}^{#1}}}
	{\ifthenelse{\isempty{#1}}
		{\int_{#2}}
		{\int_{#2}^{#1}}}
}

\NewDocumentCommand{\prodt}{mo}{%
	\mathchoice%
	{\textstyle \prod_{#1}\IfValueT{#2}{^{#2}} \displaystyle}%
	{\prod_{#1}\IfValueT{#2}{^{#2}}}%
	{\prod_{#1}\IfValueT{#2}{^{#2}}}%
	{\prod_{#1}\IfValueT{#2}{^{#2}}}%
}

\NewDocumentCommand{\prodd}{smo}{%
	\prod_{#2}\IfBooleanT{#1}{^\star}\IfValueT{#3}{^{#3}}%
}

\newcommand{\toinf}[1]{\ensuremath{#1\to\infty}\xspace}
\newcommand{\asinf}[1]{\text{as \(#1\to\infty\)}\xspace}

\newcommand{\tozero}[1]{\ensuremath{#1\to0}\xspace}

\DeclareMathOperator{\Pois}{Pois}

\DeclareMathOperator{\Bin}{Bin}

\DeclareMathOperator{\Tr}{Tr}

\newcommand{\RV}{\textsf{RV}\xspace}
\newcommand{\RVs}{\textsf{RV}s\xspace}
\newcommand{\RW}{\ifmmode \mathsf{RW} \else \textsf{RW}\xspace \fi}
\newcommand{\RWs}{\textsf{RW}s\xspace}

\newcommand{\TV}{\ifmmode \mathsf{TV} \else \textsf{TV}\xspace \fi}

\newcommand{\id}{\mathsf{id}}

\newcommand{\tmix}{t_\mix}

\newcommand{\mix}{\textnormal{mix}}

\newcommand{\MT}{\textnormal{MT}}
\newcommand{\ET}{\textnormal{ET}}

\newcommand{\Ninn}{{N\in\mathbb{N}}}

\newcommand{\mbc}{\mathbb{C}}

\newcommand{\mbn}{\mathbb{N}}

\newcommand{\mbr}{\mathbb{R}}

\newcommand{\mcl}{\mathcal{L}}
\newcommand{\mcm}{\mathcal{M}}

\newcommand{\mfG}{\mathfrak{G}}

\newcommand{\mfK}{\mathfrak{K}}

\newcommand{\mfX}{\mathfrak{X}}

\newcommand{\nt}{\addtocounter{equation}{1}\tag{\theequation}}

\newcommand{\eps}{\varepsilon}

\usepackage{manyfoot}
\SetFootnoteHook{\hspace*{-1.8em}}
\DeclareNewFootnote{bl}[gobble]
\newcommand{\blfootnote}[1]{\footnotebl{\sffamily#1}}

\def\IfAmpersandUseAlign#1#2&#3\EndIfAmpersandUseAlign%
	{%
		\if\relax\detokenize{#3}\relax
		\begin{equation*}%
		#1%
		\end{equation*}%
		\else
		\begin{align*}%
		#1%
		\end{align*}%
		\fi
	}
	\def\[#1\]%
	{%
		\IfAmpersandUseAlign{#1}#1&\EndIfAmpersandUseAlign
	}

\usepackage{eqparbox}

\usepackage[capitalise]{cleveref}

\crefname{figure}{Figure}{Figures}

\numberwithin{equation}{section}

\crefformat{equation}{\textup{(#2#1#3)}}
\crefmultiformat{equation}{\textup{(#2#1#3}}{\textup{, #2#1#3)}}{\textup{, #2#1#3}}{\textup{, #2#1#3)}}
\crefrangeformat{equation}{\textup{(#3#1#4--#5#2#6)}}

\crefformat{section}{\S#2#1#3}
\crefformat{subsection}{\S#2#1#3}
\crefformat{subsubsection}{\S#2#1#3}

\newenvironment{Proof}[1][\proofname]{%
	\proof[\upshape\bfseries\sffamily\boldmath#1]
}{\endproof}

\usepackage{etoolbox}
\usepackage{needspace}

\newcommand{\nextresult}{%
	\setcounter{introthm}{\value{introthm}}
	\setcounter{introlem}{\value{introthm}}
	\setcounter{introdefn}{\value{introthm}}
	\setcounter{intrormkT}{\value{introthm}}
}

\newtheoremstyle{sfsl}
	{1\baselineskip}		
	{1\baselineskip}		
	{\slshape}				
	{}						
	{\bfseries\sffamily}	
	{.}						
	{0.5em}					
	{\thmname{#1}\thmnumber{ #2}\thmnote{ {\mdseries(#3)}}}

\newtheoremstyle{sfup}
	{1\baselineskip}		
	{1\baselineskip}		
	{\upshape}				
	{}						
	{\bfseries\sffamily}	
	{.}						
	{0.5em}					
	{\thmname{#1}\thmnumber{ #2}\thmnote{ {\mdseries(#3)}}}

\theoremstyle{sfsl}

\newtheorem*{thm*}{Theorem}
\newtheorem{thm} {Theorem}[section]
\crefname{thm}{Theorem}{Theorems}

\newtheorem*{introthm*}{Theorem}
\newtheorem{introthm}{Theorem}

\crefname{introthm}{Theorem}{Theorems}

\newtheorem*{cor*}{Corollary}
\newtheorem{cor} [thm]{Corollary}
\crefname{cor}{Corollary}{Corollaries}

\newtheorem*{introcor*}{Corollary}

\crefname{introcor}{Corollary}{Corollaries}

\newtheorem*{introconj*}{Conjecture}
\newtheorem{introconj}{Conjecture}

\crefname{introconj}{Conjecture}{Conjectures}

\newtheorem*{introques*}{Question}

\crefname{introques}{Question}{Questions}

\newtheorem*{lem*}    {Lemma}
\newtheorem{lem} [thm]{Lemma}
\crefname{lem}{Lemma}{Lemmas}

\newtheorem*{introlem*}{Lemma}
\newtheorem{introlem}{Lemma}

\crefname{introlem}{Lemma}{Lemmas}

\newtheorem*{prop*}    {Proposition}
\newtheorem{prop} [thm]{Proposition}
\crefname{prop}{Proposition}{Propositions}

\newtheorem*{clm*}    {Claim}

\crefname{clm}{Claim}{Claims}

\newtheorem*{defn*}    {Definition}
\newtheorem{defn} [thm]{Definition}
\crefname{defn}{Definition}{Definitions}

\newtheorem*{introdefn*}{Definition}
\newtheorem{introdefn}{Definition}

\crefname{introdefn}{Definition}{Definitions}

\newtheorem*{alg*}{Algorithm}

\crefname{alg}{Algorithm}{Algorithms}

\crefname{nota}{Notation}{Notations}

\newtheorem{innercustomgeneric}{\customgenericname}
\providecommand{\customgenericname}{}

\newcommand{\customgenname}{}

\newtheorem*{conj*}   {Conjecture}

\crefname{conj}{Conjecture}{Conjectures}

\newenvironment{conj-ind*}
	{\begin{quote}\textsf{\textbf{Conjecture.}}\slshape}
	{\end{quote}}
\newenvironment{conj-ind}
	{\begin{quote}\vspace{-\glueexpr\baselineskip+\topsep}\begin{customconj}}
	{\end{customconj}\end{quote}}

\newenvironment{question-ind*}
	{\begin{quote}\textsf{\textbf{Question.}}\slshape}
	{\end{quote}}
\newenvironment{question-ind}
	{\begin{quote}\vspace{-\glueexpr\baselineskip+\topsep}\begin{customquestion}}
	{\end{customquestion}\end{quote}}

\newenvironment{openproblem-ind*}
	{\begin{quote}\textsf{\textbf{Open Problem.}}\slshape}
	{\end{quote}}
\newenvironment{openproblem-ind}
	{\begin{quote}\vspace{-\glueexpr\baselineskip+\topsep}\begin{customopenproblem}}
	{\end{customopenproblem}\end{quote}}

\newtheorem*{hypothesis*}{Hypothesis}

\newtheorem*{hyp*}{Hypothesis}

\crefname{hyp}{Hypothesis}{Hypotheses}

\newtheorem*{rmk*}{Remark}

\theoremstyle{sfup}

\crefname{defn} {Definition}{Definitions}
\crefname{defnT}{Definition}{Definitions}

\newtheorem*{defnTT}{Definition}
\newenvironment{defnt*}
	{\pushQED{\qed}\defnTT}
	{\popQED\enddefnTT}

\newtheorem{rmkT}[thm]{Remark}
	\newenvironment{rmkt}
	{\pushQED{\qed}\rmkT}
	{\popQED\endrmkT}
\crefname{rmk} {Remark}{Remarks}
\crefname{rmkT}{Remark}{Remarks}

\newtheorem*{rmkTT}{Remark}
\newenvironment{rmkt*}
	{\pushQED{\qed}\rmkTT}
	{\popQED\endrmkTT}

\crefname{rmks} {Remarks}{Remarks}
\crefname{rmksT}{Remarks}{Remarks}

\newtheorem*{rmks*} {Remarks}
\newtheorem*{rmksTT}{Remarks}
\newenvironment{rmkst*}
	{\pushQED{\qed}\rmksTT}
	{\popQED\endrmksTT}

\newtheorem{intrormkT}{Remark}
	\newenvironment{intrormkt}
	{\pushQED{\qed}\intrormkT}
	{\popQED\endintrormkT}

\crefname{intrormk} {Remark}{Remarks}
\crefname{intrormkT}{Remark}{Remarks}

\newtheorem*{intrormk*} {Remark}
\newtheorem*{intrormkTT}{Remark}
\newenvironment{intrormkt*}
	{\pushQED{\qed}\intrormkTT}
	{\popQED\endintrormkTT}

\crefname{exm} {Example}{Examples}
\crefname{exmT}{Example}{Examples}

\newtheorem*{exm*} {Example}
\newtheorem*{exmTT}{Lemma}
	\newenvironment{exmt*}
	{\pushQED{\qed}\exmTT}
	{\popQED\endexmTT}

\newtheorem*{note*} {Note}
\newtheorem*{noteTT}{Note}
	\newenvironment{notet*}
	{\pushQED{\qed}\noteTT}
	{\popQED\endnoteTT}


\newcounter{parentnumber}

\makeatletter
\newenvironment{subtheorem-num}[1]{%
	\def\subtheoremcounter{#1}%
	\refstepcounter{#1}%
	\protected@edef\theparentnumber{\csname the#1\endcsname}%
	\setcounter{parentnumber}{\value{#1}}%
	\setcounter{#1}{0}%
	\expandafter\def\csname the#1\endcsname{\theparentnumber.\arabic{#1}}%
	\expandafter\def\csname theH#1\endcsname{thm.\theparentnumber.\arabic{#1}}%
	\unskip\ignorespaces
}{%
	\setcounter{\subtheoremcounter}{\value{parentnumber}}%
	\ignorespacesafterend
}
\makeatother

\makeatletter
\let\save@mathaccent\mathaccent
\newcommand*\if@single[3]{%
  \setbox0\hbox{${\mathaccent"0362{#1}}^H$}%
  \setbox2\hbox{${\mathaccent"0362{\kern0pt#1}}^H$}%
  \ifdim\ht0=\ht2 #3\else #2\fi
  }
\newcommand*\rel@kern[1]{\kern#1\dimexpr\macc@kerna}
\newcommand*\widebar[1]{\@ifnextchar^{{\wide@bar{#1}{0}}}{\wide@bar{#1}{1}}}
\newcommand*\wide@bar[2]{\if@single{#1}{\wide@bar@{#1}{#2}{1}}{\wide@bar@{#1}{#2}{2}}}
\newcommand*\wide@bar@[3]{%
  \begingroup
  \def\mathaccent##1##2{%
    \let\mathaccent\save@mathaccent
    \if#32 \let\macc@nucleus\first@char \fi
    \setbox\z@\hbox{$\macc@style{\macc@nucleus}_{}$}%
    \setbox\tw@\hbox{$\macc@style{\macc@nucleus}{}_{}$}%
    \dimen@\wd\tw@
    \advance\dimen@-\wd\z@
    \divide\dimen@ 3
    \@tempdima\wd\tw@
    \advance\@tempdima-\scriptspace
    \divide\@tempdima 10
    \advance\dimen@-\@tempdima
    \ifdim\dimen@>\z@ \dimen@0pt\fi
    \rel@kern{0.6}\kern-\dimen@
    \if#31
      \overline{\rel@kern{-0.6}\kern\dimen@\macc@nucleus\rel@kern{0.4}\kern\dimen@}%
      \advance\dimen@0.4\dimexpr\macc@kerna
      \let\final@kern#2%
      \ifdim\dimen@<\z@ \let\final@kern1\fi
      \if\final@kern1 \kern-\dimen@\fi
    \else
      \overline{\rel@kern{-0.6}\kern\dimen@#1}%
    \fi
  }%
  \macc@depth\@ne
  \let\math@bgroup\@empty \let\math@egroup\macc@set@skewchar
  \mathsurround\z@ \frozen@everymath{\mathgroup\macc@group\relax}%
  \macc@set@skewchar\relax
  \let\mathaccentV\macc@nested@a
  \if#31
    \macc@nested@a\relax111{#1}%
  \else
    \def\gobble@till@marker##1\endmarker{}%
    \futurelet\first@char\gobble@till@marker#1\endmarker
    \ifcat\noexpand\first@char A\else
      \def\first@char{}%
    \fi
    \macc@nested@a\relax111{\first@char}%
  \fi
  \endgroup
}
\makeatother

\usepackage{xparse}
\ExplSyntaxOn
\NewDocumentCommand{\mref}{m}{\quinn_mref:n {#1}}
\seq_new:N \l_quinn_mref_seq
\cs_new:Npn \quinn_mref:n #1
{
	\seq_set_split:Nnn \l_quinn_mref_seq { , } { #1 }
	\seq_pop_right:NN \l_quinn_mref_seq \l_tmpa_tl
	( 
	\seq_map_inline:Nn \l_quinn_mref_seq
	{ \ref{##1},\nobreakspace } 
	\exp_args:NV \ref \l_tmpa_tl 
	) 
}
\ExplSyntaxOff

\usepackage[doi=false,isbn=false,url=false,
	hyperref=auto,
	sorting=nyt,
	maxnames=10,
	maxcitenames=3,
	backend=biber,
	texencoding=auto,
	giveninits=true,
	block=space,
]{biblatex}

\DefineBibliographyStrings{english}{%
	andothers = {et al}
} 
 
\DeclareFieldFormat{eprint:arxiv}{%
	\em \href{https://arxiv.org/abs/#1}{arXiv:\allowbreak #1}}

\DeclareFieldFormat{eprint:mrnumber}{%
	\href{http://www.ams.org/mathscinet-getitem?mr=MR#1}{MR\allowbreak #1}}

\DeclareFieldFormat{eprint:jstor}{%
	\href{http://www.jstor.org/stable/#1}{JSTOR\allowbreak #1}}

\DeclareFieldFormat{eprint:customeprint}{%
	\em Available at \upshape \ttfamily \href{https://#1}{#1}}

\DeclareFieldFormat{eprint:inprep}{%
	\em In preparation}

\DeclareFieldFormat{eprint:note}{%
	\em {#1}}

\DeclareFieldFormat{eprint:onarxiv}{%
	\em Available on arXiv}

\DeclareFieldFormat{eprint:toappear}{%
	\mbox{\em To appear}}

\DeclareSourcemap{
	\maps[datatype=bibtex]{
		\map{
			\step[fieldsource=mrnumber,		fieldtarget=eprint, final]
			\step[fieldset=eprinttype,		fieldvalue=mrnumber]
		}
		\map{
			\step[fieldsource=arxiv,		fieldtarget=eprint, final]
			\step[fieldset=eprinttype,		fieldvalue=arxiv]
		}
		\map{
			\step[fieldsource=jstor,		fieldtarget=eprint, final]
			\step[fieldset=eprinttype,		fieldvalue=jstor]
		}
		\map{
			\step[fieldsource=customeprint,	fieldtarget=eprint, final]
			\step[fieldset=eprinttype,		fieldvalue=customeprint]
		}
		\map{
			\step[fieldsource=inprep,		fieldtarget=eprint, final]
			\step[fieldset=eprinttype,		fieldvalue=inprep]
		}
		\map{
			\step[fieldsource=note,			fieldtarget=eprint, final]
			\step[fieldset=eprinttype,		fieldvalue=note]
		}
		\map{
			\step[fieldsource=onarxiv,		fieldtarget=eprint, final]
			\step[fieldset=eprinttype,		fieldvalue=onarxiv]
		}
		\map{
			\step[fieldsource=toappear,		fieldtarget=eprint, final]
			\step[fieldset=eprinttype,		fieldvalue=toappear]
		}
	}
}

\DeclareFieldFormat{doi}{%
	\href{https://doi.org/#1}{DOI}%
}

\DeclareFieldFormat{url}{%
	\href{https://#1}{\texttt{#1}}%
}

\DeclareFieldFormat{comments}{%
	\em #1%
}

\DeclareFieldFormat
	[article,inbook,incollection,inproceedings,patent,thesis,unpublished]
	{title}{{#1\isdot}}

\DeclareFieldFormat
	[article,book,inbook,incollection,inproceedings,patent,thesis,unpublished, online]
	{date}{\mkbibparens{#1}}
\DeclareFieldFormat
	[article,book,inbook,incollection,inproceedings,patent,thesis,unpublished]
	{volume}{\mkbibbold{#1}}
\DeclareFieldFormat
	{pages}{\mkbibparens{#1}}

\DeclareFieldFormat{edition}{%
	\ifinteger{#1}
	{\mkbibordedition{#1}~\bibstring{edition}\addcomma}
	{#1~\bibstring{edition}\addcomma}}

\renewbibmacro*{volume+number+eid}{%
	\printfield{volume}%
\setunit*{\adddot}%
	\printfield{number}%
\setunit{\space}%
	\printfield{eid}%
}

\DeclareBibliographyDriver{article}{%
	\usebibmacro{bibindex}%
	\usebibmacro{begentry}%
	\usebibmacro{author}%
\setunit{\addspace}%
	\usebibmacro{date}%
\newunit\newblock
	\usebibmacro{title}%
\newunit\newblock
	\printfield{journaltitle}\addperiod%
\setunit*{\addspace}%
	\usebibmacro{volume+number+eid}%
\setunit*{\addspace}\newblock
	\printfield{pages}
\setunit*{\addspace}
	\printfield{comments}%
	\usebibmacro{eprint}%
\setunit*{\addnbspace}%
	\printfield{doi}%
	\usebibmacro{finentry}%
}

\DeclareBibliographyDriver{online}{%
	\usebibmacro{bibindex}%
	\usebibmacro{begentry}%
	\usebibmacro{author}%
\setunit{\addspace}%
	\usebibmacro{date}%
\newunit\newblock
	\usebibmacro{title}%
\newunit\newblock
	\usebibmacro{eprint}%
	\usebibmacro{finentry}%
}

\DeclareBibliographyDriver{book}{%
	\usebibmacro{bibindex}%
	\usebibmacro{begentry}%
	\usebibmacro{author}%
\setunit{\addspace}\newblock
	\usebibmacro{date}%
\newunit\newblock
	\usebibmacro{title}%
\setunit*{\addspace}\newblock
	\printfield{edition}%
\newunit\newblock
	\printlist{publisher}
\setunit*{\addspace}%
	\printfield{volume}%
\setunit*{\addspace}\newblock%
	\printfield{comments}%
	\usebibmacro{eprint}%
\setunit*{\addspace}
	\printfield{doi}
	\usebibmacro{finentry}%
}

\DeclareBibliographyDriver{thesis}{%
	\usebibmacro{bibindex}%
	\usebibmacro{begentry}%
	\usebibmacro{author}%
\setunit{\addspace}\newblock
	\usebibmacro{date}%
\newunit\newblock
	\usebibmacro{title}.%
\setunit*{\addspace}\newblock
	\space Thesis, \printlist{institution}
	\usebibmacro{eprint}%
	\printfield{comments}%
\setunit*{\addspace}
	\printfield{doi}%
	\usebibmacro{finentry}%
}

\DeclareBibliographyDriver{inproceedings}{%
	\usebibmacro{bibindex}%
	\usebibmacro{begentry}%
	\usebibmacro{author}%
\setunit{\addspace}%
	\usebibmacro{date}%
\newunit\newblock
	\usebibmacro{title}%
\newunit\newblock
	\printfield{booktitle}\addcomma%
\setunit*{\addspace}%
	\printfield{series}\addcomma%
\setunit*{\addspace}\newblock
	\usebibmacro{editor}\addcomma
\setunit*{\addspace}\newblock
	\printlist{publisher}
\setunit*{\addspace}%
	\printfield{volume}%
\setunit*{\addspace}%
	\printfield{pages}
\setunit*{\addspace}
	\usebibmacro{eprint}%
	\printfield{comments}%
\setunit*{\addnnspace}
	\printfield{doi}
	\usebibmacro{finentry}%
}

\DeclareBibliographyDriver{inbook}{%
	\usebibmacro{bibindex}%
	\usebibmacro{begentry}%
	\usebibmacro{author}%
\setunit{\addspace}%
	\usebibmacro{date}%
\newunit\newblock
	\usebibmacro{title}%
\newunit\newblock
	\printfield{booktitle}\addcomma%
\setunit*{\addspace}%
	\printfield{series}\addcomma%
\setunit*{\addspace}\newblock
	\usebibmacro{editor}\addcomma
\setunit*{\addspace}\newblock
	\printlist{publisher}
\setunit*{\addspace}%
	\printfield{volume}%
\setunit*{\addspace}%
	\printfield{pages}
\setunit*{\addspace}
	\usebibmacro{eprint}%
	\printfield{comments}%
\setunit*{\addspace}
	\printfield{doi}
	\usebibmacro{finentry}%
}

\DeclareBibliographyDriver{incollection}{%
	\usebibmacro{bibindex}%
	\usebibmacro{begentry}%
	\usebibmacro{author}%
\setunit{\addspace}%
	\usebibmacro{date}%
\newunit\newblock
	\usebibmacro{title}%
\newunit\newblock
	\printfield{booktitle}\addcomma%
\setunit*{\addspace}%
	\printfield{series}\addcomma%
\setunit*{\addspace}\newblock
	\printlist{publisher}
\setunit*{\addspace}%
	\printfield{volume}%
\setunit*{\addspace}%
	\printfield{pages}
\setunit*{\addspace}
	\usebibmacro{eprint}%
	\printfield{comments}%
\setunit*{\addspace}
	\printfield{doi}
	\usebibmacro{finentry}%
}

\AtEveryBibitem{%
	\clearfield{day}%
	\clearfield{month}%
} 
 
\DeclareNameAlias{sortname}{given-family}
\DeclareNameAlias{default}{given-family}

\newcommand{\cdf}{\textsf{cdf}\xspace}
\newcommand{\mgf}{\textsf{mgf}\xspace}
\newcommand{\mgfs}{\textsf{mgf}s\xspace}

\newcommand{\whp}{\textsf{whp}\xspace}

\newcommand{\DS}{Diaconis--Shahshahani\xspace}

\DeclareMathOperator{\Fix}{Fix}
\newcommand{\symgr}{\mathfrak S}

\newcommand{\trans}{\mathfrak T}

\DeclareMathOperator{\HG}{HG}

\newcommand{\kk}{k}

\renewcommand{\AA}{A}
\renewcommand{\aa}{a}
\newcommand{\bb}{b}
\newcommand{\BB}{B}
\newcommand{\cc}{c}

\newcommand{\ee}{\mathrm{e}}

\newcommand{\ff}{f}
\newcommand{\FF}{F}
\newcommand{\GG}{G}
\newcommand{\hh}{h}
\newcommand{\ii}{i}
\newcommand{\II}{I}
\newcommand{\jj}{j}
\newcommand{\JJ}{J}
\newcommand{\KK}{K}
\newcommand{\LL}{L}
\newcommand{\mm}{m}
\newcommand{\MM}{M}
\newcommand{\nn}{n}
\newcommand{\NN}{N}

\newcommand{\PP}{P}

\newcommand{\rr}{r}

\renewcommand{\SS}{S}
\renewcommand{\tt}{t}
\newcommand{\TT}{T}

\newcommand{\UU}{U}
\newcommand{\vv}{v}
\newcommand{\VV}{V}
\newcommand{\ww}{w}
\newcommand{\WW}{W}
\newcommand{\xx}{x}
\newcommand{\XX}{X}
\newcommand{\yy}{y}
\newcommand{\YY}{Y}
\newcommand{\zz}{z}
\newcommand{\ZZ}{Z}

\newcommand{\CEP}{\ifmmode \mathsf{CEP} \else \textsf{CEP}\xspace \fi}
\newcommand{\EP}{\ifmmode \mathsf{EP} \else \textsf{EP}\xspace \fi}
\newcommand{\IP}{\ifmmode \mathsf{IP} \else \textsf{IP}\xspace \fi}
\newcommand{\RT}{\ifmmode \mathsf{RT} \else \textsf{RT}\xspace \fi}

\newcommand{\unif}[1][]{%
	\mathcal U%
	\ifthenelse{\equal{#1}{}}{}{_{#1}}%
}
\NewDocumentCommand{\Unif}{sm}{%
	\IfBooleanTF{#1}%
		{\mathcal U_{#2}}%
		{\operatorname{Unif}(#2)}%
}

\newcommand{\ft}{\widehat}
\newcommand{\proj}{\widetilde}
\newcommand{\lift}{}

\NewDocumentCommand{\dtv}{smm}{%
	d_\TV%
	\IfBooleanT{#1}{\bigl}(#2, \: #3\IfBooleanT{#1}{\bigr})%
}

\NewDocumentCommand{\dd}{sm}{%
	d_{#2\IfBooleanT{#1}{^\star}}%
}

\NewDocumentCommand{\DD}{smooo}{%
	d_{#2}%
	\IfValueT{#3}{%
		\IfBooleanT{#1}{\bigl}(%
		#3%
		\IfValueT{#4}{, \: #4}%
		\IfBooleanT{#1}{\bigr})%
	}%
}

\NewDocumentCommand{\lam}{s}{%
	\lambda\IfBooleanT{#1}{^\star}%
}

\NewDocumentCommand{\qq}{m}{%
	q_{#1}%
}

\usepackage{empheq}

\title{\sffamily%
	Limit Profile for Projections of Random Walks on Groups%
}

\author{\sffamily Evita Nestoridi\quad Sam Olesker-Taylor}
\date{}

\begin{document}

\maketitle

\blfootnote{%
	Evita Nestoridi%
\quad
	\href{mailto:evrydiki.nestoridi@stonybrook.edu}{evrydiki.nestoridi@stonybrook.edu}%
\hfill%
	\href{mailto:sam.olesker-taylor@warwick.ac.uk}{sam.olesker-taylor@warwick.ac.uk}%
\quad%
	Sam Olesker-Taylor%
\\
	Department of Mathematics, Stony Brook University%
\hfill
	Department of Statistics, University of Warwick%
\\
	Supported by NSF Grant DMS-2052659%
\hfill
	Partially supported by EPSRC Grant EP/N004566/1%
}

\vspace{-6ex}

\renewcommand{\abstractname}{\sffamily Abstract}

\begin{abstract}\noindent
Establishing cutoff—an abrupt transition from `not mixed' to `well mixed'—is a classical topic in the theory of mixing times for Markov chains. Interest has grown recently in determining not only the \emph{existence} of cutoff and the order of its mixing time and window, but the exact \emph{shape}, or \emph{profile}, of the convergence inside the window. Classical techniques, such as coupling or $\ell_2$-bounds, are typically too crude to establish this and there has been a push to develop general techniques \cite{T:limit-profile,NOt:limit-profiles:rev,N:limit-profiles:comp}.

We build upon this work, extending from conjugacy-invariant random walks on groups to certain projections. We exemplify our method by analysing the $k$-particle interchange process on the complete $n$-graph with $k \asymp n$. This is a projection of the random-transposition card shuffle, which corresponds to $k = n$, analysed in \cite{T:limit-profile}.
\end{abstract}

\small
\begin{quote}
\begin{description}
	\item [Keywords:]
	cutoff,
	limit profile,
	spectral theory,
	random walk on groups,
	projections of random walks,
	representation theory,
	Fourier transform,
	homogeneous space
	
	\item [MSC 2020 subject classifications:]
	20C15, 20C30; 43A30, 43A65; 60B15, 60C05, 60J10
	
\end{description}
\end{quote}
\normalsize






\sffamily
\boldmath

\setcounter{tocdepth}{1}

\makeatletter
\renewcommand\tableofcontents{%
	\@starttoc{toc}%
}
\makeatother
\section*{Table of Contents}
\vspace*{-2ex}

\tableofcontents

\unboldmath
\normalfont




\section{Introduction}

Developing a theory for determining the limit profiles of Markov chains is a very recent topic in studying Markov chains.
The limit profiles of certain popular Markov chains have been determined a while ago: eg,
	random walk on hypercube \cite{DGB:hypercube-profile,NOt:limit-profiles:rev},
	riffle shuffles \cite{BD:riffle-shuffle},
	exclusion process on the circle \cite{L:exclusion-circle-profile}
and
	Ramanujan graphs \cite{LP:ramanujan}.
Some more recent examples include
	random transpositions and $k$-cycles~\cite{T:limit-profile,NOt:limit-profiles:rev},
	star transpositions \cite{N:limit-profiles:comp},
	biased card shuffles~\cite{Z:cutoff-biased-shuffle},
	quantum random transpositions~\cite{FTW:quantum-rt},
	some random walks on random Cayley graphs~\cite{HOt:rcg:abe:cutoff},
	the asymmetric exclusion process~\cite{BN:asep-profile},
	urn models \cite{NOt:limit-profiles:rev},
	repeated averages~\cite{CDSZ:repeated-avg}
and
	a Moran model~\cite{C:moran-limit-profile}.
All these references are somewhat isolated, case-by-case proofs, without an underlying technique or theory.
In many cases, cutoff was known well before the profile.
This is partially an artefact of the diversity of profile functions:
	some are
		Poissonian~\cite{T:limit-profile} or Gaussian~\cite{HOt:rcg:abe:cutoff};
	others involve more complicated distributions, such as
		Tracy--Widom~\cite{Z:cutoff-biased-shuffle} or free Meixner~\cite{FTW:quantum-rt}.

Rectifying this lack of general methodology, spectral techniques, which can be applied to a variety of chains, have been introduced recently \cite{T:limit-profile,NOt:limit-profiles:rev,N:limit-profiles:comp}.
\textcite{T:limit-profile} kicked things off, introducing an approximation technique, applying it to prove that the total-variation distance to uniformity of the random-transposition shuffle is
\[
	\DD*\TV[\Pois(1+\ee^{-\cc})][\Pois(1)]
\Quad{after}
	\tfrac12 \nn \rbr{ \log \nn + \cc }
\quad
	\text{steps asymptotically},
\]
by extending Fourier-analytic techniques of \textcite{DS:random-trans} for conjugacy-invariant random walks (\textit{\RW}s).
We extended this from transpositions ($2$-cycles) to $\kk$-cycles \cite{NOt:limit-profiles:rev}.
We extended the approximation technique beyond conjugacy-invariant \RWs to so-called \textit{Gelfand pairs} and applied it to the multi-urn Ehrenfest diffusion model \cite{NOt:limit-profiles:rev}.
We also extended to reversible Markov chains and applied it to a statistical-physics Gibbs sampler.
A comparison method was introduced in
\cite{N:limit-profiles:comp},
	allowing the limit profile to be deduced from a sufficiently similar chain,
and applied it to star transpositions.

The purpose of this paper is to set up a technique for finding the limit profile of a Markov chain which can be viewed as a projection of another.
Gelfand pairs mentioned above are a special case of this.
We apply the technique to the $\kk$-particle interchange process on the complete $\nn$-graph, which is a projection of the random-transposition shuffle.
	Roughly, it is
	\[
		\DD*\TV[\Pois\rbb{\tfrac \kk\nn + \ee^{-\cc}}][\Pois\rbb{\tfrac \kk\nn}]
	\Quad{after}
		\tfrac12 \nn \rbr{ \log \kk + \cc }
	\quad
		\text{steps asymptotically if}
	\quad
		\tfrac \kk\nn \asymp 1.
	\]

\subsection{$k$-Particle Interchange Process and Mixing Times}

We study the limit profile of the \textit{$\kk$-particle interchange process} (\textit{$\kk$-\IP}).

\begin{introdefn}[$\kk$-\IP]
Place $\kk \le \nn$ distinct particles on the vertices an $\nn$-graph $\GG$.
Each vertex gets at most one particle and therefore there are always $\nn - \kk$ empty sites.
A single step of the \textit{$\kk$-particle interchange process} $\IP(\kk, \GG)$ dynamics chooses an edge uniformly at random:
\begin{itemize}[noitemsep]
	\item 
	if both endpoints have a particle,
	then swap the positions of the two particles;
	
	\item 
	if only one endpoint does,
	then move that particle to the other endpoint;
	
	\item 
	if neither endpoint has a particle,
	then do nothing.
\end{itemize}

When $\GG = K_\nn$ is the complete $\nn$-graph,
we abbreviate $\IP(\kk, K_\nn)$ as $\kk$-\IP.
Our complete graphs always have self-loops.
Nothing happens upon selection of a self-loop:
	it is a lazy~step.
\end{introdefn}

We write $\alpha \cq \kk/\nn$ for the proportion of sites which are occupied.
Our main result evaluates the limit profile for the $\kk$-\IP on the complete $\nn$-graph $K_\nn$ when $\alpha \asymp 1$:
	it is
	\[
		\DD*\TV[\Pois\rbr{\alpha + \ee^{-\cc}}][\Pois(\alpha)]
	\Quad{after}
		\tfrac12 \nn \rbr{ \log \kk + \cc }
	\quad
		\text{steps}.
	\]

The $\kk$-\IP is a natural projection of the random-transposition shuffle:
	the latter equivalent to $\nn$-\IP
and
	$\kk$-\IP is obtained from the $\nn$-\IP by following only $\kk$ cards.

The objective of this paper is to introduce a method for analysing a \RW on a space $\XX$ which can be viewed as a projection of some simpler \RW on a group $\GG$.
In particular, we study \RWs on \textit{homogeneous spaces} $\XX = \GG/\KK$, where $\KK \le \GG$ is a subgroup of $\GG$.
We derive theory for general $(\GG, \KK)$.
For $\kk$-\IP,
	$\GG = \symgr_\nn$
and
	$\KK = \symgr_{\nn-\kk} \times \symgr_1 \times ... \times \symgr_1$.


\medskip

We carefully define \textit{mixing}, \textit{cutoff} and \textit{limit profiles} now. Let $\Omega$ be a finite set and $\PP$ a transition matrix on $\Omega$.
Then, $\PP^\tt(\xx,\yy)$ is the probability of moving from $\xx$ to $\yy$ in $\tt$ steps,
for all $\xx,\yy \in \Omega$ and all $\tt \in \mbn$.
If $\PP$ is irreducible and aperiodic, then the basic limit theorem of Markov chains
(eg, \cite[Theorem~1.8.3]{N:markov})
tells us that $\PP^\tt(\xx,\cdot)$ converges to the (unique) invariant distribution $\pi$ as \toinf \tt with respect to the \textit{total variation} (\textit{\TV}) distance $\DD\TV$:
\[
	\DD*\TV[\PP^\tt(\xx,\cdot)][\pi]
\cq
	\tfrac12 \sumt{\yy \in \Omega} \abs{P^\tt(\xx,\yy) - \pi(\yy)}
\Qfor
	\xx \in \Omega
\Qand
	\tt \in \mbn;
\]
ie, \TV is half the $\ell_1$ distance.
The (\textit{worst-case}) $\eps$-\textit{mixing time} is then defined,
	\mbox{for $\eps \in [0,1]$,~%
by}
\[
	\tmix(\eps)
\cq
	\inf\brb{\tt \ge 0 \midb d_\TV(\tt) \le \eps}
\Qwhere
	\DD\TV[\cdot]
\cq
	\maxt{\xx \in \Omega}
	\DD\TV[\PP^\tt(\xx,\cdot)][\pi].
\]

Take a sequence of Markov chains indexed by $\NN$,
and denote the \TV distance for the $\NN$-th chain by $\DD\TV^{(N)}(\cdot)$.
If there exists $(\tt_\star^{(\NN)}, \ww_\star^{(\NN)})_\Ninn$ such that
$\lim_{\NN\to\infty} \ww_\star^{(\NN)} / \tt_\star^{(\NN)} = 0$,
\[
	\LIM{\alpha \to -\infty}
	\LIM{\toinf \NN} \,
	\DD\TV^{(\NN)}\rbb{ \tt_\star^{(\NN)} + \alpha \ww_\star^{(\NN)}}
=
	1
\Qand
	\LIM{\alpha \to +\infty}
	\LIMSUP{\toinf \NN} \,
	\DD\TV^{(\NN)}\rbb{ \tt_\star^{(\NN)} + \alpha \ww_\star^{(\NN)} }
=
	0,
\]
then the sequence exhibits \textit{cutoff} at
\(
	(\tt_\star^{(\NN)})_{\NN\in\mbn}
\)
with \textit{window} order at most
\(
	(\ww_*^{(\NN)})_{\NN\in\mbn}.
\)

The existence of cutoff is a fascinating question.
The \textit{product condition}, requiring the product of the \textit{spectral gap} and \textit{mixing time} to diverge, is a necessary condition.
It was proved to be sufficient in $\ell_p$-distance mixing with $p > 1$ by \textcite{CSc:cutoff-ergodic}.
It is not sufficient in \TV, ie $\ell_1$, in complete generality---see, eg, \cite[Example~18.7]{LPW:markov-mixing}---but is conjectured so be sufficient in many natural scenarios.
This leaves open the question of characterising the chains for which the product condition implies cutoff in \TV.
The implication has been established for birth-and-death chains \cite{DLP:cutoff-birth-death} and, more generally, random walks on trees~\cite{BHP:cutoff-rev}, as well as exclusion processes with reservoirs~\cite{S:cutoff-ep-reservoir}.

Significant progress for general chains has been made by \citeauthor{S:cutoff-nonneg-curv} recently.
They gave a sufficient condition for cutoff, based on the concept of \textit{varentropy}, and a condition for \textit{non-negatively curved} Markov chains to satisfy this condition, in \cite{S:cutoff-nonneg-curv}.
Very recently, they showed that this varentropy condition is in fact sharp on sparse expanders in \cite{S:varentropy-sharp}.

One can look beyond the cutoff time and window to determine the \textit{limit profile}:
	\[
		\varphi(\alpha)
	\cq
		\LIM{N \to \infty} \,
		\DD\TV^{(N)}\rbb{ \tt_\star^{(N)} + \alpha w_\star^{(N)} }
	\Qfor
		\alpha \in \mbr,
	\]
if the limit exists.
There has been much less investigation into the limit profile than cutoff.

\subsection{Limit Profile for $k$-\IP on the Complete $n$-Graph}

Cutoff for $k$-\IP on $K_n$ was established by Lacoin and Leblond~\cite[Theorem~1.2]{LL:cutoff-ep-ip} whenever $1 \ll k \ll \sqrt n$, with the correct window.
Their lower bound is valid for all $k$ which,
combined with the well-known upper bound for $k = n$, establishes cutoff with the correct window whenever $k \asymp n$.
However, this approach is not refined enough to obtain the limit profile.

\begin{introthm}[Poisson Profile]
\label{res:intro:res:ip}
Let $\alpha \in (0,1]$.
Suppose that $(\nn_\NN)_{\NN\in\mbn}$ and $(\kk_\NN)_{\NN\in\mbn}$ are diverging sequences of integers with $1 \le \kk_\NN \le \nn_\NN$ for all $\NN$ and $\kk_\NN / \nn_\NN \to \alpha$ as $\NN \to \infty$.~%
Set
\[
	\tt_\cc(\nn, \kk)
\cq
	\tfrac12 \nn \rbr{ \log \kk + \cc }
=
	\tfrac12 \nn \rbr{ \log \nn + \cc + \log(k/n) }
\Qforall
	\cc \in \mbr.
\]
Then, the \TV distance $\DD\TV^{(\NN)}(\cdot)$ to equilibrium for $\kk_\NN$-\IP on $\KK_{\nn_\NN}$ at time $\tt_\cc(\nn_\NN, \kk_\NN)$ satisfies
\[
	\DD*\TV^{(\NN)}\rbb{ \tt_\cc(\nn_\NN, \kk_\NN) }
\to
	\DD*\TV[\Pois(\alpha + \ee^{-\cc})][\Pois(\alpha)]
\quad
	\asinf \NN.
\]
	%
\end{introthm}

\begin{intrormkt}
As with many of these processes based on random transpositions, the limit profile is the same as that observed by looking only at the number of fixed points---here, amongst the $k = \alpha n$ particles.
Specifically, the probability that a particle is in its original~place~is
\[
	\rbb{ 1 - (\nn-1)/\binomt \nn2}^{\tt_\cc}
+	\tfrac1\nn
+	\oh{\tfrac1\nn}
=
	\tfrac1\nn \rbb{ 1 + \ee^{-\cc}/\alpha } + \oh{\tfrac1\nn}.
\]
The total number amongst the $k = \alpha \nn$ labelled particles in the original place is then
\[
	\Pois\rbb{ \alpha \nn \rbb{ \tfrac1\nn( 1 + \ee^{-\cc}/\alpha ) + \oh{\tfrac1\nn} } }
\to
	\Pois\rbr{ \alpha + \ee^{-\cc} }.
\]
On the other hand, in a uniform permutation, it is $\Pois(\alpha)$.
This all assumes that $\alpha \asymp 1$.
\end{intrormkt}

We believe that the profile is Gumbel when $1 \ll \kk \ll \nn$.
This is a natural extension since
\[
	\DD*\TV[\Pois(e^{-c})][\Pois(0)]
=
	\exp{-e^c},
\]
which is the Gumbel tail, under the convention $\Pois(0) = 0$.

\refstepcounter{introconj}
\begin{introconj}[Gumbel Profile]
Suppose that $(k_N)_{N\in\mbn}$ and $(n_N)_{N\in\mbn}$ are diverging sequences of integers with $1 \le \kk_\NN \le \nn_\NN$ for all $\NN$ and $k_N \to \infty$ and $k_N/n_N \to 0$ as $N \to \infty$.~%
Set
\[
	t_c(n,k)
\cq
	\tfrac12 n ( \log k + c )
\Qfor
	\cc \in \mbr.
\]
Then, the \TV distance $\DD\TV^{(\NN)}(\cdot)$ to equilibrium for $\kk_\NN$-\IP on $\KK_{\nn_\NN}$ at time $\tt_\cc(\nn_\NN, \kk_\NN)$ satisfies
\[
	\DD*\TV^{(\NN)}\rbb{ \tt_\cc(\nn_\NN, \kk_\NN) }
\to
	\exp{ - e^{-\cc} }
\quad
	\asinf \NN.
\]
\end{introconj}

\begin{intrormkt}
Our argument should extend to \emph{some} $\alpha \ll 1$.
A number of the approximations become more technically challenging in this set-up and it is not clear how far the computations can be pushed.
More details are given in \cref{rmk:ip:main:k-small} in \cref{sec:ip:main:just} after the calculations.

If $k$ is sufficiently small%
---specifically, if $k^2 \log k \ll n$---%
then all $k$ particles will jump to an (almost) uniform space independently before interacting with each other, with high probability.
The mixing time is thus the time to touch all the particles.
\begin{itemize}[noitemsep, topsep = \smallskipamount]
	\item 
	If $k \gg 1$, then this coupon-collector time concentrates, and has Gumbel fluctuations.
	
	\item 
	If $k \asymp 1$, then there are finitely many particles, and there is no concentration or cutoff.
\end{itemize}

In summary, the fluctuations should be Gumbel when $1 \ll k \ll \sqrt{n/\log n}$ or when $k$ is sufficiently close to $n$.
We conjecture that they are Gumbel whenever $1 \ll k \ll n$.
\end{intrormkt}

The key trick in finding this limit profile is to couple $\kk$-\IP with $\nn$-\IP with random initial state.
From a representation-theoretic point of view, this allows us to delay evaluating certain complicated Fourier transforms for $\kk$-\IP, rather working with the simpler ones for $\nn$-\IP.

Obtaining only cutoff is possible working with $\kk$-\IP directly; see \cite{S:many-urn-bl}.
Indeed, the traces are calculable as the Fourier transforms are still diagonal, but no longer multiples of the identity.
The profile requires more refined information on the diagonals than just~the~sum.
This extra information seems unmanageable.
Hence the need for the $\kk$-\IP-to-$\nn$-\IP simplification.

\begin{intrormkt*}
	%
The case $\alpha = 1$, ie $\nn$-\IP, is the random-transposition shuffle, analysed recently by \textcite{T:limit-profile}.
The case $\alpha \to 1$ is an easy deduction from the fixed-$\alpha$ case.
In this case,
\[
	\tt_\cc(\nn, \nn) - \tt_\cc(\nn, \kk)
=
	\tfrac12 \nn \log \alpha
\ll
	\nn
\]
and $\nn$ is the order of the cutoff window.
\TV is monotone under projection and the $\kk'$-\IP can be obtained from the $\kk$-\IP, with $\kk' \le \kk$, by forgetting about $\kk-\kk' \ge 0$ of the particles.
The claim follows by sandwiching the $\kk$-\IP between $\kk'$-\IP and $\nn$-\IP where $\kk' \cq \alpha' \nn$ for some fixed $\alpha' < 1$, arbitrarily close to $1$.
Herein, we assume that $\alpha$ is a fixed element of $(0, 1)$.
\end{intrormkt*}

\subsection{Approximation Lemma for \RWs on Homogeneous Spaces}
\label{sec:intro:approx}

Cutoff and the $\ell_2$ limit-profile for the $\kk$-\IP on $K_\nn$ was established by \textcite{FJ:spectrum-ep-ip} by computing the entire spectrum.
The $\ell_2$ distance requires only the eigenvalues, not the eigenfunctions.
Establishing the \TV limit-profile requires the eigenfunctions.
\textcite[Lemma~2.1]{T:limit-profile} introduced an approximation lemma for doing precisely this in the set-up on conjugacy-invariant \RWs.
One of the main contributions of the current article is an extension of \citeauthor{T:limit-profile}'s approximation lemma; see \cref{res:intro:approx}.
We need to introduce some standard representation-theoretic terminology in order to set-up this lemma.

\nextresult

\begin{subtheorem-num}{introdefn}
	\label{def:intro:approx}

\begin{subtheorem-num}{intrormkT}
	\label{rmk:intro:approx}

First and foremost, we introduce the most fundamental part of representation theory.

\begin{introdefn}[Representations]
	Let $\GG$ be a finite group and $\VV$ a finite-dimensional vector space over $\mbc$.
	A \textit{representation} $\rho$ of $\GG$ over $\VV$ is an action
	\(
		(g, \vv) \mapsto \rho(g) \cdot \vv
	:
		\GG \times \VV \to \VV
	\)
	such that $\rho(g) : \VV \to \VV$ is an invertible, linear map for all $g \in \GG$.
	Let
		$\dd\rho \cq \dim \VV$ denote the \textit{dimension}
	and
		$\chi_\rho(\cdot) \cq \Tr\rbr{ \rho(\cdot) }$ the \textit{character}
	of $\rho$.
	A representation $(\rho, \VV)$~is~\textit{irreducible}~if~its only $\GG$-invariant vector subspaces are trivial:
		$\WW \le \VV$\! and $\rho(G) \WW \le \WW$\! implies $\WW \in \bra{ \bra{0}, \VV }$.
	Let
		$\ft \GG$ be the set of irreducible representations (\textit{irreps})
	and	
		$\ft \GG^\star$ the set of non-trivial irreps.
\end{introdefn}

A \textit{left coset} of a subgroup $K$ of a group $G$ is of the form $g K = \bra{ gk \mid k \in K }$, for $g \in G$.
Left-cosets $g K$ and $h K$ are \textit{equivalent} if $g^{-1} h \in K$.
A set of \textit{coset-representatives} is a collection $\bra{ g_1 K, ..., g_N K }$ of pairwise non-equivalent cosets, where $N = \abs{G/K} = \abs G / \abs K$.

We are going to project \RWs from a group $G$ to a \textit{homogeneous space} of left costs.


\begin{introdefn}[Homogeneous Space]
\label{def:intro:approx:hom}
	Let $\GG$ be a group and $\KK \le \GG$ a subgroup.
	The \textit{homogeneous space} $\XX \cq \GG / \KK$ is the set of left-cosets of $\GG$ mod $\KK$.
	We associate $\XX$ with an arbitrary set $\bra{\xx_1, ..., \xx_\NN} \subseteq \GG$ of coset representatives for $\XX$, including the identity $\id \in G$.
	
	Write $\mcm(\XX)$ for the set of probability measures on~$\XX$.
	For $\proj \nu, \proj{\nu_0} \in \mcm(\XX)$,
	let $\RW(\proj \nu, \proj{\nu_0})$ denote the \RW on $\XX$ driven by $\proj \nu$ and started from $\proj{\nu_0}$:
		its time-$\tt$ law is $\proj{\nu_0} * \proj \nu^{*t}$ for all $\tt \in \mbn$.
	
\end{introdefn}

We require a way of lifting measures from the homogeneous space to the group.

\begin{introdefn}[Lifting and Conjugacy Invariance]
	Let $\proj \nu \in \mcm(\XX)$.
	Define $\lift \nu \in \mcm(\GG)$ by lifting $\proj \nu$ to $\GG$:
	\(
		\lift \nu(g)
	\cq
		\proj \nu(\xx) / \abs K
	\)
	where $\xx$ is the unique coset-representative such that $g = \xx \kk$ for some $\kk \in \KK$.
	In particular, $\delta_\id \in \mcm(\XX)$ is lifted to $\Unif{\KK} \in \mcm(\GG)$.
	Say that $\proj \nu$ is \textit{$\KK$-conjugacy invariant} if $\lift \nu$ is conjugacy invariant---ie, constant on each conjugacy class of $\GG$.
	In this case, for a representation $\rho$ of $\GG$, define the \textit{character ratio}
	\(
		\qq\rho
	\cq
		\dd\rho^{-1}
		\sumt{g \in \GG}
		\lift \nu(g) \chi_\rho(g).
	\)
\end{introdefn}

We now state the main approximation lemma, which is used for finding the limit profile.

\begin{introlem}[Approximation Lemma]
\label{res:intro:approx}
	Let $\XX = \GG/\KK$ be a homogeneous space and $\SS \subseteq \ft \GG^\star$.
	Let $\proj \nu \in \mcm(\XX)$ be $\KK$-conjugacy invariant.
	Let $\UU \sim \Unif\XX$ and $\VV \sim \Unif\KK$.
	Let $\tt \ge 0$.~%
	Then,
	\[
		\absb{ 
			\DD*\TV[\RW_\tt(\proj \nu, \delta_\id)][\Unif \XX]
		-	\tfrac12
			\ex{ 
				\absb{ 
					\ex{ 
						\sumt{\rho \in \SS}
						\dd \rho
						\qq\rho^\tt
						\chi_\rho(\UU \VV)
					\mid
						\UU
					}
				}
			}
		}
	\le
		\tfrac12
		\sumt{\rho \in \GG^\star \setminus \SS}
		\sqrt{ \mm_\rho \dd\rho }
		\abs{\qq\rho^\tt},
	\]
	where
		$\mm_\rho$ is the multiplicity of the trivial representation inside the restriction of $\rho$ to $\KK$
	and
		the subscript-$\tt$ in $\RW_\tt(\cdot, \cdot\cdot)$ indicates that the \RW is run for time $\tt$.
\end{introlem}

\begin{intrormkt}
\label{rmk:intro:approx:tv-pres}
	The \TV distance $\DD\TV[\nu][\pi]$ can be realised by an optimal coupling.
	From this, it is immediate that lifting and projecting both preserve the \TV distance:
	\[
		\DD*\TV[\RW_\tt(\proj \nu, \proj{\nu_0})][\Unif\XX]
	=
		\DD*\TV[\RW_\tt(\lift \nu, \lift{\nu_0})][\Unif\GG]
	\Qforall
		\tt \ge 0.
	\]
	It thus does not matter whether we view a process as a \RW on the homogeneous space $\XX$ or lifted to the group $\GG$:
		both have \emph{exactly} the same \TV distance to equilibrium.
\end{intrormkt}

\begin{intrormkt}
	Taking $\KK \cq \bra{\id}$ in \cref{res:intro:approx} recovers \cite[Lemma~2.1]{T:limit-profile}.
	We obtained an extension in the set-up of \textit{Gelfand pairs} in \cite[Lemma~C]{NOt:limit-profiles:rev}, not requiring conjugacy invariance:
		the characterising property of a Gelfand pair is that $\mm_\rho \in \bra{0, 1}$ for all $\rho \in \ft \GG$.
\end{intrormkt}

\end{subtheorem-num}

\end{subtheorem-num}

\section{Random Walks and Fourier Analysis}
\label{sec:rw-ft}

Here,
	$\GG$ denotes a finite group,
	$\ft \GG$ the set of irreducible representations (\textit{irreps}) of $\GG$
and
	$\ft \GG$ the non-trivial irreps.
	$\KK \le \GG$ is a subgroup of $\GG$
and
	$\XX \cq \GG / \KK$ the corresponding homogeneous space.
Given a set $\YY$, the set of probability measures on $\YY$ is denoted $\mcm(\YY)$.
%
We use the $\ell_1$, rather than \TV, distance, to avoid carrying factors of $\tfrac12$:
	$\DD1[\cdot][\cdot\cdot] \cq 2 \DD\TV[\cdot][\cdot\cdot]$.

The uniform distribution is always invariant for \RWs on groups.
We abbreviate $\Unif*\YY \cq \Unif\YY$
and drop it from the relative distances:
	if $\nu \in \mcm(Y)$, then write $\DD{\cdot}[\nu] \cq \DD{\cdot}[\nu][\Unif*\YY]$.

\subsection{Random Walks on Groups and Homogeneous Spaces}
\label{sec:rw-ft:rw}

The fundamental tool in representation theory of finite groups is the \textit{Fourier transform}.
%
	\[
		\ft \nu(\rho)
	\cq
		\sumt{g \in \GG}
		\nu(g) \rho(g)
	\Qfor
		\rho \in \ft \GG.
	\]

\begin{prop}[Fourier Inversion Formula; {\cite[\S 3.10]{CsST:harmonic-analysis-finite-groups}}]
\label{res:rw-ft:rw:fourier-inversion}
	Let $\nu \in \mcm(\GG)$.
	Then,
	\[
		\DD1[\nu][\Unif*\GG]
	=
		\abs \GG^{-1}
		\sumt{g \in \GG}
		\absb{ 
				\sumt{\rho \in \ft \GG^\star}
				\dd\rho
				\Tr\rbb{ \rho(g) \ft \nu(\rho) }
			}.
	\]
\end{prop}

The Fourier inversion formula is particularly amenable to \RWs on groups,
as the Fourier transform turns convolutions into products.
This means that \cref{res:rw-ft:rw:fourier-inversion} gives
\[
	\DD*1[\RW_\tt(\nu, \nu_0)]
=
	\abs \GG^{-1}
	\sumt{g}
	\absb{
		\sumt*{\rho}
		\dd\rho \Tr\rbb{ \rho(g) \ft \nu(\rho)^\tt \ft{\nu_0}(\rho) }
	},
\]
where $\sumt*{\rho} \cq \sumt{\rho \in \GG^\star}$ and $\sumt{g} \cq \sumt{g \in \GG}$.
If, further, $\nu$ is constant on conjugacy classes, a standard application of Schur's lemma gives $\ft \nu(\rho) = \qq\rho I$ for some $\qq\rho \in \mbc$, where $I$ is the $\dd\rho \times \dd\rho$ identity matrix.
This can then be pulled out of the trace as just a constant factor:
\[
	\DD*1[\RW_\tt(\nu, \nu_0)]
=
	\abs \GG^{-1}
	\sumt{g}
	\absb{
		\sumt*{\rho}
		\dd\rho \qq\rho^\tt \Tr\rbb{ \rho(g) \ft{\nu_0}(\rho) }
	}.
\]
If $\nu_0$ is a point-mass measure, which is the usual set-up, then $\ft{\nu_0}$ is the identity map.
Then,
\[
	\DD*1[\RW_\tt(\nu, \nu_0)]
=
	\abs \GG^{-1}
	\sumt{g}
	\absb{
		\sumt*{\rho}
		\dd\rho \qq\rho^\tt \chi_\rho(g)
	}
=
	\exb{ \absb{ \sumt*{\rho} \dd\rho \qq\rho^\tt \chi_\rho(\UU) } }
\Qwhere
	\UU \sim \Unif*\GG.
\]

The set-ups studied in this paper do not fall into this ``conjugacy-invariant \RW started from a point-mass'' regime, however.
Rather, they are driven by a projection of such a measure from $\GG$ to $\GG / \KK$.
Recall that this projection is \TV-preserving from \cref{rmk:intro:approx:tv-pres}:
\[
	\DD*\TV[\RW_\tt(\proj \nu, \proj{\nu_0})][\Unif\XX]
=
	\DD*\TV[\RW_\tt(\lift \nu, \lift{\nu_0})][\Unif\GG].
\]
It thus suffices to study
	$\RW(\lift{\nu_t}, \lift{\nu_0})$ instead of $\RW(\proj{\nu_t}, \proj{\nu_0})$.
The former is, in our set-ups, a conjugacy-invariant \RW on $\GG$ started from the non-deterministic measure $\lift{\nu_0}$.
Now,
\[
	\nu_t = \nu_0 * \nu^{*t},
\quad
	\ft{\nu}(\rho) = \qq\rho I
\Qand
	\proj{\nu_0} = \delta_{\id_\XX},
\Quad{so}
	\ft{\nu_t}(\rho)
=
	\qq\rho \ft{\nu_0}(\rho)
\Qand
	\nu_0 = \Unif*\KK.
\]

The key to our method is in avoiding \emph{actually evaluating} the Fourier transform
\[
	\ft{\nu_0}(\rho)
=
	\ft{\Unif*\KK}(\rho)
=
	\abs \KK^{-1}
	\sumt{\kk \in \KK}
	\rho(\kk).
\]
Rather, we just use this expanded-summation form.
The following lemma is now immediate.

\begin{lem}
\label{res:rw-ft:exact}
	Let $\XX = \GG / \KK$ be a homogeneous space.
	Let $\proj \nu \in \mcm(\XX)$ be $\KK$-conjugacy invariant with character ratios $\rbr{\qq\rho \mid \rho \in \ft \GG}$.
	Let $\proj{\nu_0} \cq \delta_{\id_\XX} \in \mcm(\XX)$.
	Let $\tt \ge 0$.
	Then,
	\[
	&	\DD*1[\RW_\tt(\proj \nu, \proj{\nu_0})][\Unif*\XX]
	=
		\DD*1[\RW_\tt(\lift \nu, \lift{\nu_0})][\Unif*\GG]
	\\&\qquad
	=
		\abs \GG^{-1}
		\sumt{g \in \GG}
		\absb{ 
			\abs \KK^{-1}
			\sumt{\kk \in \KK}
			\sumt{\rho \in \GG^\star}
			\dd\rho \qq\rho^\tt
			\chi_\rho(g \kk)
		}
	\\&\qquad
	=
		\ex{ 
			\absb{ 
				\ex{ 
					\sumt*\rho
					\dd\rho \qq\rho^\tt
					\chi_\rho(\UU \VV)
				\mid
					\UU
				}
			}
		}
	\Qwhere
		(\UU, \VV) \sim \Unif*\GG \times \Unif*\KK,
	\]
	It is equivalent to let $\UU \sim \Unif*\XX$, viewing $\XX = \bra{\xx_1, ..., \xx_\NN} \subseteq \GG$.
\end{lem}

The motivation for this lemma was evaluation of the limit profile for random transpositions given by \textcite{T:limit-profile}.
There is no projection in their set-up; this is equivalent to $\KK = \bra{\id}$---no averaging over $\KK$.
They approximately evaluate the inner sum
\(
	\sumt*\rho
	\dd\rho \qq\rho^\tt
	\chi_\rho(\sigma)
\)
for (almost) any $\sigma \in \symgr_\nn$, the symmetric group on $\nn$ objects.
They then take an absolute value and average over $\sigma \in \symgr_\nn$.
We must average over $\KK$ before taking the absolute value.


\subsection{The Approximation Lemma}
\label{sec:rw-ft:approx}

The \DS upper bound \cite{DS:random-trans} is a standard tool for upper bounding the distance to uniformity for a conjugacy-invariant \RW:
\(
	\DD1[\nu_\tt]^2
\le
	\sumt*\rho
	\dd\rho^2 \qq\rho^{2\tt}.
\)
\textcite{T:limit-profile} improved on this by first separating out the sum over irreps into `main' and `error' terms.

We use the same idea.
The main challenge in the proof of \cref{res:intro:approx} is in the error.~%
Simple applications of the triangle inequality and Cauchy--Schwarz suffice for conjugacy-invariant \RWs and Gelfand pairs---ie, when $K = \bra{\id}$ or $m_\rho \in \bra{0,1}$---but not in our general set-up.


\begin{Proof}[Proof of Lemma \ref{res:intro:approx}]
We start from the $\ell_1$ representation of \cref{res:rw-ft:exact}:
\[
&	\absb{
		\DD*\TV[\RW_\tt(\proj \nu, \delta_\id)][\Unif \XX]
	-	\exb{
			\absb{
				\exb{
					\sumt{\rho \in \SS}
					\dd\rho \qq\rho^\tt \chi_\rho(\UU \VV)
					\mid
					\UU
				}
			}
		}
	}
\\&\qquad
=
	\absb{
		\exb{
			\absb{
				\exb{
					\sumt*{\rho}
					\dd\rho \qq\rho^\tt \chi_\rho(\UU \VV)
					\mid
					\UU
				}
			}
		-	\absb{
				\exb{
					\sumt{\rho \in \SS}
					\dd\rho \qq\rho^\tt \chi_\rho(\UU \VV)
					\mid
					\UU
				}
			}
		}
	}
\\&\qquad
\le
	\absb{
		\exb{
			\absb{
				\sumt*{\rho \notin \SS}
				\dd\rho \qq\rho^\tt
				\exb{ \chi_\rho(\UU \VV) \mid \UU }
			}
		}
	}
\le
	\sumt*{\rho \notin \SS}
	\dd\rho \abs{\qq\rho^\tt}
	\exb{ \absb{ \ex{ \chi_\rho(\UU \VV) \mid \UU } } },
\]
where $(\UU, \VV) \sim \Unif*\GG \times \Unif*\KK$, as before.
We return to the original $\ft{\Unif*\KK}(\rho)$ viewpoint,
rather than the expanded $\abs \KK^{-1} \sumt{\kk \in \KK} \rho(\kk)$ viewpoint,
to bound the error:
\[
	\ex{ \chi_\rho(g \VV) }
=
	\Tr\rbb{
		\rho(g)
	\cdot
		\abs \KK^{-1}
		\sumt{\kk \in \KK}
		\rho(\kk)
	}
=
	\Tr\rbb{
		\rho(g)
		\ft{\Unif*\KK}(\rho)
	}.
\]

Let $(\rho, \WW)$ be an irrep of $G$ and let $\WW^\KK$ denote the space of $\KK$-invariant vectors in $\WW$:
\[
	\WW^\KK
\cq
	\bra{ \ww \in \WW \mid \rho(\kk) \ww = \ww \text{ for all } \kk \in \KK }.
\]
Let $\mm_\rho$ denote its dimension:
\(
	\mm_\rho
\cq
	\dim \WW^\KK.
\)
Choose an orthonormal basis $\bra{\ww_1, ..., \ww_{\dd\rho}}$ of $\WW$
such that the first $\mm_\rho$ vectors are $\KK$-invariant, ie $\bra{\ww_1, ..., \ww_{\mm_\rho}} \subseteq \WW^\KK$.
Thus, $\bra{\ww_1, ..., \ww_{\mm_\rho}}$ forms an orthonormal basis for $\WW^\KK$.
Then, in this basis, the Fourier transform
\[
	\ft{\Unif*\KK}(\rho)
=
	\abs \KK^{-1}
	\sumt{\kk \in \KK}
	\rho(\kk)
\Quad{is represented as}
{\large
\begin{psmallmatrix}
	I & 0 \\
	0 & 0
\end{psmallmatrix}}
\in
	\mbc^{\dd\rho \times \dd\rho}	
\]
where $I \in \mbc^{\mm_\rho \times \mm_\rho}$ is the $\mm_\rho \times \mm_\rho$ identity matrix.
Define the \textit{matrix coefficients}
\[
	\rho_{\ii,\jj}(g)
\cq
	\langle \rho(g) \ww_\jj, \: \ww_\ii \rangle
\Qfor
	\ii, \jj \in [\dd\rho]
\Qand
	g \in \GG;
\]
here, $\langle \cdot, \: \cdot\cdot \rangle$ is the usual inner product between vectors.
With this notation,
\[
	\Tr\rbb{ \rho(g) \ft{\Unif*\KK}(\rho) }
=
	\sumt{\ii \in [\mm_\rho]}
	\rho_{\ii,\ii}(g)
\Qforall
	g \in \GG,
\]
recalling that the trace is basis-independent.
The matrix coefficients are orthogonal:
\[
	\ex{ \rho_{\ii,\ii'}(\UU) \overline{\rho_{\jj,\jj'}(\UU)} }
=
	\abs \GG^{-1}
	\langle \rho_{\ii,\ii'}, \: \rho_{\jj,\jj'} \rangle
=
	\delta_{\ii,\jj} \delta_{\ii',\jj'} / \dd\rho
\Qforall
	\ii,\ii',\jj,\jj' \in [\dd\rho];
\]
see \cite[Lemma~3.6.3]{CsST:harmonic-analysis-finite-groups}.
We use this orthogonality after applying Cauchy--Schwarz:
\[
&	\ex{ \absb{ \sumt{\ii \in [\mm_\rho]} \rho_{\ii,\ii}(\UU) } }{}^2
\le
	\ex{ \absb{ \sumt{\ii \in [\mm_\rho]} \rho_{\ii,\ii}(\UU) }^2 }
\\&\qquad
=
	\sumt{\ii \in [\mm_\rho]}
	\ex{ \rho_{\ii,\ii}(\UU) \overline{\rho_{\ii,\ii}(\UU)} }
+	\sumt{\ii, \jj \in [\mm_\rho] : \ii \ne \jj}
	\ex{ \rho_{\ii,\ii}(\UU) \overline{\rho_{\jj,\jj}(\UU)} }
=
	\mm_\rho / \dd\rho.
\]

The claim now follows immediately by combining the upper bounds just established.
\end{Proof}


\begin{rmkt}
	\textcite{S:many-urn-bl} observed that the \DS upper bound gives
	\[
		\DD*\TV[\RW_\tt(\proj \nu, \delta_\id)][\Unif \XX]{}^2
	\le
		\tfrac14
		\sumt*{\rho}
		\mm_\rho \dd\rho \abs{\qq\rho^{2\tt}}.
	\qedhere
	\]
\end{rmkt}

%
%

\section{Preliminaries for the $\kk$-Particle Interchange Process}
\label{sec:prelim}

\subsection{$\kk$-\IP as \RW on Homogeneous Space}
\label{sec:prelim:hom}

We have now proved the general approximation lemma for ``projected conjugacy-invariant \RWs''.
We now turn to our application,
	which is that of the $\kk$-\IP on the complete $\nn$-graph.
The following lemma sets up the more general \textit{coloured exclusion process} as such a \RW.

\begin{defn}[Coloured Exclusion Process]
	Let $\nn, \NN \in \mbn$ and $\aa = (\aa_1, ..., \aa_\NN) \in \mbn^\NN$ with $\sumt{\jj\ge1} \aa_\jj = \nn$; let $\aa_0 \cq 0$.
	Fix an $\nn$-vertex graph $\GG$ and place $\nn$ particles, labelled $1$ through $\nn$, disjointly in its vertices.
	Paint particles labelled $\sumt{\jj' \le \jj-1} \aa_{\jj'} + 1$ through $\sumt{\jj' \le \jj} \aa_{\jj'}$ with colour $\jj$, for each $\jj \in [\NN]$.
%
	A single step of the \textit{coloured exclusion process} $\CEP(\aa, \GG)$ chooses an edge  uniformly at random:
		the particles at either endpoint of the edge are swapped.
\end{defn}

Let $\trans_\nn \subseteq \symgr_\nn$ denote the set of all transpositions in $\symgr_\nn$.
This is a conjugacy class.

\begin{lem}
	Let $\nn, \NN \in \mbn$ and $\aa = (\aa_1, ..., \aa_\NN) \in \mbn^\NN$ with $\sumt{\jj\ge1} \aa_\jj = \nn$.
	Let $\mfG \cq \symgr_\nn$ and
	\(
		\mfK
	\cq
		\symgr_{\aa_1} \times \cdots \times \symgr_{\aa_\NN};
	\)
	let $\mfX \cq \mfG / \mfK$.
	Let $\lift \mu_\RT \cq \Unif{\trans_\nn}$ and $\lift \mu \cq \tfrac{\nn-1}{\nn} \mu_\RT + \tfrac1\nn \delta_\id \in \mcm(\mfG)$;
	let $\proj \mu \in \mcm(\mfX)$ denote the projection of $\mu$ from $\mfG$ to $\mfX$.
	Then,
	\(
		\CEP(\aa, K_\nn)
	=
		\RW\rbr{ \proj \mu, \delta_\id }
	\)
\end{lem}

\begin{Proof}
	%
Particles of the same colour are indistinguishable. Thus, in the lifted space $\GG = \symgr_\nn$, we may permute the $\aa_\jj$ particles of colour $\jj$ arbitrarily and not affect the distribution.
	%
\end{Proof}


The interchange process with $\kk$ particles is simply the coloured exclusion process where $\kk$ particles are given a unique, distinguishing colour and the remainder the same colour.

\begin{lem}
	Let $\GG$ be an $\nn$-graph and $\kk \le \nn$.
	Then, $\IP(\kk, \GG) = \CEP\rbr{ (\nn-\kk, 1, ..., 1), \GG }$.
\end{lem}


\subsection{Representation Theory for the Symmetric Group}
\label{sec:prelim:rep}

The $\kk$-\IP is fundamentally a projection of the random-transposition shuffle, which is a \RW on the symmetric group $\symgr_\nn$.
Representation theory for $\GG = \symgr_\nn$ is well understood.

The set $\ft \GG = \ft{\symgr_\nn}$ of irreps is indexed by partitions of $\nn$.
We abuse notation~slightly by writing $\lambda \in \ft{\symgr_\nn}$ to denote the irrep canonically associated to the partition $\lambda \vdash \nn$ of $\nn$.

The set $\trans_\nn \subseteq \symgr_\nn$ of transpositions is a conjugacy class, so the Fourier transform $\ft{\mu_\RT}(\lambda)$ is a multiple of the identity, by Schur's lemma.
Taking traces, this multiple is $\rr_\lambda \cq \chi_\lambda(\tau) / \dd\lambda$, where $\tau \in \trans_\nn$ is any transposition.
A simple-to-evaluate expression for $\rr_\lambda$ is known.

\begin{lem}[{\cite[Theorem~10.6.2]{CsST:harmonic-analysis-finite-groups}}]
\label{res:prelim:rep:ft-rt}
	Let $\lambda \vdash \nn$.
	Write $\lambda' \vdash \nn$ for the \textit{transpose} of $\lambda$.
	Then,
	\[
	\textstyle
		\rr_\lambda
	=
		\frac1{\nn(\nn-1)}
		\sumd{\ii \in [\nn]}
		\lambda_\ii \rbb{ \lambda_\ii - (2\ii - 1) }
	=
		\binom \nn2^{-1}
		\sumd{\ii \in [\nn]}
		\rbb{ \binom{\lambda_\ii}2 - \binom{\lambda'_\ii}2 }.
	\]
\end{lem}

This is the fundamental building block for calculating the Fourier transform of our driving measures.
We also need some bounds on the dimensions of the irreducible representations.

\begin{lem}
\label{res:prelim:rep:d_lambda}
	Let $\lambda \vdash \nn$.
	For a box $(\ii, \jj)$ in the Young diagram of $\lambda$, let
	\(
		\hh_{(\ii,\jj)}
	\cq
		\rbr{ \lambda_\ii  - j }
	+	\rbr{ \lambda'_\jj - i }
	+	1
	\)
	denote the \textit{hook length} of the box $(\ii, \jj)$.
	Then, the \textit{hook-length formula} says that
	\[
		\dd\lambda
	=
		n!
	\big/
		\prodt{(\ii,\jj) \in \lambda}
		\hh_{(\ii,\jj)}.
	\]
	In particular, writing $\rr \cq \nn - \lambda_1$, we have
	\[
		\binomt \nn\rr
		\dd*\lambda
		(1 - 2\rr/\nn)
	\cdot
		\one{\rr \le \tfrac12 \nn}
	\le
		\dd\lambda
	\le
		\binomt \nn\rr
		\dd*\lambda.
	\]
\end{lem}


\begin{Proof}
The hook-length formula is a standard result; see, eg, \cite[Theorem~4.2.14]{CsST:harmonic-analysis-finite-groups}.

The bounds follow from writing the quotient $\dd\lambda / \dd*\lambda$ using the hook-length formula:
\[
	\frac{\dd\lambda}{\dd*\lambda}
=
	\frac{(\nn)_{\lambda_1}}{\prodt{\ii=1}[\lambda_1] (\lambda_1 - \ii + \lambda'_\ii)}
=
	\binom \nn{\lambda_1}
	\prodd{\ii=1}[\lambda_1]
	\frac{\lambda_1 - \ii + \lambda'_\ii}{\lambda_1 - \ii + 1}.
\]
Clearly, $\lambda'_\ii \ge 1$ for all $\ii \le \lambda_1$.
Thus, $\dd\lambda \le \binom \nn\rr \dd*\lambda$.
On the other hand, to minimise the product, each term should be as close to equal as possible.
Hence, given $\rr = \nn-\lambda_1 \le \tfrac12 \nn$, it is minimised by $\lambda = (\nn-\rr, \rr)$:
then, $\lam* = (\rr)$ and $\dd*\lam = 1$;
so,
\[
	\dd\lambda
\ge
	\dd{(\nn-\rr, \rr)}
=
	\binom \nn\rr
	\frac{\nn-2\rr+1}{\nn-\rr+1}
=
	\binom \nn\rr
	\rbbb{ 1 - \frac \rr{\nn-\rr+1} }
\ge
	\binom \nn\rr
	\rbbb{ 1 - \frac{2\rr}{\nn} }.
\]

The upper bound can also be proved combinatorially.
First, choose which $\rr = \abs{\lam*}$ numbers are placed in $\lam*$; then, place the remainder in the first row, in increasing order.
This gives an upper bound since not all enumerations will give rise to increasing columns.
\end{Proof}

Calculating, or estimating, the multiplicities $\mm_\lambda$ is fundamentally important for bounding the error term. We recall that $\mm_\lambda$ is the multiplicity of the irrep $\lambda$ in the restriction of $\lambda$ from $\GG = \symgr_\nn$ to $\KK = \symgr_{\aa_1} \times \cdots \times \symgr_{\aa_\NN}$ in the general \CEP set-up;
	$\NN = \kk+1$, $\aa_1 = \cdots = \aa_\kk = 1$ and $\aa_{\kk+1} = \nn-k$ for the $\kk$-\IP.
In particular, the zero-multiplicity irreps can be ignored.

\begin{defn}[Young Tableaux, {\cite[\S 3.7.1]{CsST:rep-theory-Sn}}]
	Let $\lambda, \mu \vdash \nn$.
	A \textit{tableau of shape $\lambda$ and type $\mu$} is a filling of the Young diagram associated to the partition $\lambda$ with integers:
		the integer $\jj$ appears $\mu_\jj$ times for each $\jj$.
	It is \textit{semi-standard} if the integers form
		weakly increasing sequences along the rows
	and
		strictly increasing sequences along the columns.
	It is \textit{standard} if the rows are strictly increasing, which means that $\mu = (1, ..., 1) \vdash \nn$.
\end{defn}

The multiplicity $\mm_\lambda$ is given by the Young rule, which we describe now.

\begin{lem}[Young Rule, {\cite[Theorem~3.7.10 and Corollary~3.7.11]{CsST:rep-theory-Sn}}]
\label{res:prelim:rep:mult-cep}
	The multiplicity $\mm_\lambda$ equals the number of semi-standard tableaux of shape $\lambda$ and type $(\aa_1, ..., \aa_\NN)$.
\end{lem}

This is easy to calculate or bound in the special cases of the $\kk$-\IP.

\begin{lem}[Multiplicity for $\kk$-\IP]
\label{res:prelim:rep:mult-ip}
	Let $\lambda \vdash \nn$ and $\rr \cq \nn-\lambda_1$.
	Let $\mm_\lambda$ denote the number of semi-standard of shape $\lambda$ and type $(\nn-\kk, 1, ..., 1) \vdash \nn$;
	always, $\mm_\lambda \le \dd\lambda$.
	Then, $\mm_\lambda > 0$ if and only if $\rr \le \kk$.
	If $\rr \le \kk$, then
	\(
		\mm_\lambda
	\le
		\binom \kk\rr \dd*\lambda;
	\)
	further, this is an equality if $\nn-\kk \ge \rr$.
\end{lem}

\begin{Proof}
We must fill the Young diagram associated to $\lambda$ with the numbers $1$ through $\kk+1$,
	with $1$ appearing $\nn-\kk$ times and the remainder precisely once.
The integers must form a weakly increasing sequence along the rows and \emph{strictly} increasing along the columns.
The value $\dd\lambda$ corresponds to placing numbers $1$ through $\nn$%
	---ie, $\kk = \nn-1$.
Trivially, $\mm_\lambda \le \dd\lambda$.

The $\nn-\kk$ $1$s must be placed in the first row of the $\lambda$-diagram;
Therefore if $\lambda_1 < \nn-\kk$, ie $\rr > \kk$, then the $1$s cannot be placed in the first row; thus $\mm_\lambda = 0$.

After we place the $\nn-\kk$ $1$s in the first row, $\rr$ of the remaining $\kk$ numbers are placed in $\lam*$. 
These $\rr$ entries in $\lam*$ must form a standard $\lam*$-tableau.
There are $\binom \kk\rr$ ways of choosing the $\rr$ out of $\kk$.
Not all of the standard $\lam*$-tableaux give a semi-standard $\lambda$-tableau, though:
	a non-$1$ entry in the first row may be larger than the entry next to it in $\lam*$.
However, every desired semi-standard tableau can be obtained in this manner.
Thus, $\mm_\lambda \le \binom \kk\rr \dd*\lam$.

The upper bound is an equality precisely when $\lambda_2 = \lam*_1 \le \nn-\kk$, as then all entries of the first row of $\lambda$ which are adjacent to $\lam*$ are filled with a $1$.
The inequality $\lambda_2 \le \nn-\kk$ is implied by $\rr \le \nn-\kk$ as $\lambda_2 \le \abs{\lam*} = \rr$.
\end{Proof}

\begin{cor}[Reweighted Multiplicity for $\kk$-\IP]
\label{res:prelim:rep:mult-reweighted}
	For all $\lambda \vdash \nn$,
	writing $\rr \cq \nn - \lambda_1$,
	we have
	\[
		\sqrt{ \mm_\lambda \dd\lambda }
	\le
		\dd*\lambda \sqrt{ \binomt \kk\rr \binomt \nn\rr }
	\le
		\dd*\lambda (\nn \kk)^{\rr/2} / \rr!.
	\]
\end{cor}

\section{Limit Profile for the Interchange Process}

This section is devoted to establishing the limit profile for $\kk$-\IP.
Recall that
\[
	\tt_\cc
\cq
	\tt_\cc(\kk, \nn)
=
	\tfrac12 \nn \rbr{ \log \kk + \cc }
=
	\tfrac12 \nn \rbr{ \log \nn + \cc + \log \alpha }
\]
is the proposed mixing time.
We start by giving brief outline of the proof now, citing the results to come.
Throughout, we work with the $\ell_1$ distance,
so as not to carry factors of $\tfrac12$.

To emphasise that we are working with the $\kk$-\IP,
we replace $(\GG, \KK, \UU, \VV)$ with $(\mfG, \mfK, \Sigma, \Pi)$:
\[
	\text{write}
\quad
	\mfG \cq \symgr_\nn
\Qand
	\mfK \cq \symgr_{\nn-\kk} \times \symgr_1 \times \cdots \times \symgr_1;
\Quad{let}
	(\Sigma, \Pi) \sim \Unif*\mfG \times \Unif*\mfK.
\]
We also use $\lambda \vdash \nn$ for the irreps, rather than $\rho$, or even $\rho_\lambda$.
Also, as the only process being studied is $\kk$-\IP, so we abbreviate its $\ell_1$ distance to uniformity after $\tt$ steps by $\DD1(\tt)$.


\begin{Proof}[Skeleton of Proof of \cref{res:intro:res:ip}]
We use \cref{res:intro:approx},
	placing partitions $\lambda$ with long first row---namely, $\rr \cq \nn - \lambda_1 \asymp 1$---in the main term
and
	the remaining in the error term:
\begin{gather*}
	\abs{
		\DD1(\tt)
	-	\MT_\tt
	}
\le
	\ET_\tt
\quad
	\text{where}
\quad
\left\{
\begin{aligned}
	\MT_\tt
&\cq
	\exb{
		\absb{
			\exb{
				\sumt*{\lambda \vdash \nn : \lambda_1 \ge \nn - \MM}
				\dd\lambda \qq\lambda^\tt \chi_\lambda(\Sigma \Pi)
			\mid
				\Sigma
			}
		}
	},
\\
	\ET_\tt
&\cq
	\sumt{\lambda \vdash \nn : \lambda_1 < \nn - \MM}
	\sqrt{ \mm_\lambda \dd\lambda }
	\abs{\qq\lambda}^\tt.
\end{aligned}
\right.
\end{gather*}
We need only consider $\lambda$ with $\mm_\lambda \ne 0$; such partitions have first row at least $\nn - \kk$ by \cref{res:prelim:rep:mult-ip}.
We then need to estimate the main term and control the error when $\tt = \tt_\cc$:
\begin{alignat*}{2}
	\ET_\tt
&=
	\oh1
&&\quad
	\text{by \cref{res:ip:error}};
\\
	\MT_\tt
&=
	\DD*1[\Pois(\alpha+\ee^{-\cc})][\Pois(\alpha)] + \oh1
&&\quad
	\text{by \cref{res:ip:main}}.
\end{alignat*}
The $\oh1$ terms here implicitly include letting \toinf \MM \asinf \nn, but arbitrarily slowly.
\end{Proof}

We make one important remark regarding notation.
At some points it is helpful to think of \emph{following} $\kk$ cards; others, \emph{randomising} $\nn-\kk$ cards.
For this reason, we introduce $\mm \cq \nn-\kk$.
We reserve $\mm$ for this value throughout.
We could use only $\nn$ and $\kk$, but a lot of the formulas later are much more natural when viewed in the second manner, thus having $\mm$ makes them more-easily interpretable.
The triple $(\nn, \kk, \mm) \in \mbn$ will \emph{always} satisfy $\nn = \mm + \kk$, so $0 \le \mm, \kk \le \nn$.
We typically do not repeat or define these in statements below.

\subsection{Character Estimation}
\label{sec:ip:char:raw}

We start by analysing the character ratio $\qq\lambda$.
We have in the back of our minds that the main term will correspond to $\lambda$ with $\rr \cq \nn-\lambda_1 \asymp 1$; we thus require accurate estimations in this regime.
We also need an upper bound valid for larger $\rr$ to control the error term.

Recall from \cref{res:prelim:rep:ft-rt}
that the character ratio $\rr_\lambda = \chi_\lambda(\tau) / \dd\lambda$ satisfies
\[
\textstyle
	\rr_\lambda
=
	\frac1{\nn(\nn-1)}
	\sumd{\ii \in [\nn]}
	\lambda_\ii \rbb{ \lambda_\ii - (2\ii - 1) }
=
	\binom \nn2^{-1}
	\sumd{\ii \in [\nn]}
	\rbb{ \binom{\lambda_\ii}2 - \binom{\lambda'_\ii}2 }.
\]
	


\begin{lem}
\label{res:ip:char:raw}
	Let $\lambda \vdash \nn$.
	The character ratio $\qq\lambda$ for the $\kk$-\IP satisfies the following relations:
	\[
		\tfrac1\nn + \tfrac{\nn-1}\nn \rr_\lambda
	=
		\qq\lambda
	&=
		\tfrac1\nn + \tfrac1{\nn^2} \sumt{\ii \in [\nn]} \lambda_\ii \rbb{ \lambda_\ii - (2\ii-1) }
	\label{eq:ip:char:raw:eq}
	\nt
	\\
		1 - \tfrac{2 \rr}{\nn} \rbb{ 1 + \tfrac{\rr}{\nn} }
	\le
		\qq\lambda
	&\le
		1 - \tfrac{\rr}{\nn}
	\label{eq:ip:char:raw:all}
	\nt
	\\
		\qq\lambda
	&\le
		1 - \tfrac{2\rr}{\nn} \rbb{ 1 - \tfrac{\rr}{\nn} }
	\Qif
		\lambda_1 \ge \tfrac12 \nn.
	\label{eq:ip:char:raw:long}
	\nt
	\]
\end{lem}

Many similar results appear in the literature;
see, e.g, \cite[Lemma~3.2]{D:group-rep}.

\begin{Proof}[Proof of \cref{res:ip:char:raw}]
There is laziness-$\tfrac1\nn$ in our driving measure corresponding to choosing the same card twice.
A non-lazy step is a random-transposition step.
Hence, \cref{eq:ip:char:raw:eq}~holds.


If $\rr = \nn$, ie $\lambda = (1^\nn)$, then $\rr_\lambda = -1$, so $\qq\lambda = -1 + \tfrac2\nn$ and the claims hold.
Next,
	\[
		\sumt{\ii \in [\nn]}
		\lambda_\ii \rbb{ \lambda_\ii - 2(\ii-1) }
	&
	=
		\sumt{\ii \in [\nn]}
		\lambda_\ii^2
	-	2 \sumt{\ii \in [\nn]}
		\lambda_\ii (\ii-1)
	\\&
	\hspace*{+0.25em}
	\mathllap{\Bigg\{}
	\hspace*{-0.25em}
	\begin{aligned}
	&	\le \nn \lambda_1 = \nn^2 - \rr \nn,
	\\
	&	\ge \lambda_1^2 - 2 \rr^2 = \nn^2 - 2 \rr \nn - \rr^2.
	\end{aligned}
	\]
	Dividing by $\binom \nn2 = \tfrac12 \nn (\nn-1)$ and using
	\(
		\tfrac{\nn-1}{\nn} + \tfrac{\rr}{\nn-1}
	\le
		1 + \tfrac{\rr}{\nn}
	\)
	when $\rr \le \nn-1$ gives \cref{eq:ip:char:raw:all} when $\lambda \ne (1^\nn) \vdash \nn$.
Finally,
	direct calculation with $\rr_\lambda$ from \cref{res:prelim:rep:ft-rt} shows that $\lambda \mapsto \rr_\lambda$ is a monotone function; see \cite[Corollary~10.6.3]{CsST:harmonic-analysis-finite-groups}.
	Hence, \cref{eq:ip:char:raw:long} holds:
	if $\rr \le \tfrac12 \nn$, then
	\[
		\qq\lambda
	\le
		\qq{(\nn-\rr, \rr)}
	=
		1 - 2 \rr (\nn-\rr+1) / \nn^2
	\le
		1 - 2 \rr (1 - \rr/\nn) / \nn.
	\qedhere
	\]
\end{Proof}

%


We need to estimate $\qq\lambda$ raised to the power of
\(
	\tt_\cc
=
	\tfrac12 \nn \rbr{ \log \kk + \cc }.
\)


\begin{lem}
\label{res:ip:char:pwr}
	Let $\cc \in \mbr$ and $\tt \cq \tt_\cc$.
	Let $\lambda \vdash \nn$ and $\rr \cq \nn - \lambda_1$.
	The following inequalities~hold:
	\begin{alignat}{2}
		\qq\lambda^\tt
	&\le
		\rbr{ \ee^{\cc} \kk }^{-\rr(1+1/200)/2}
	&\Qif
		\rr
	&\ge
		\tfrac1{100} \nn
	\Qand
		\qq\lambda \ge 0;
	\label{eq:ip:char:pwr:short}
	\nt
	\\
		\qq\lambda^\tt
	&\le
		\rbr{ \ee^{\cc} \kk }^{-\rr}
		\expb{ 2 \rr^2 \log \kk / \nn }
	&\Qif
		\rr
	&\le
		\tfrac14 \nn,
	\label{eq:ip:char:pwr:long:upper}
	\nt
	\\
		\qq\lambda^\tt
	&\ge
		\rbr{ \ee^{\cc} \kk }^{-\rr}
		\expb{ - 7 \rr^2 \log \kk / \nn }
	&\Qif
		\rr
	&\le
		\tfrac14 \nn.
	\label{eq:ip:char:pwr:long:lower}
	\nt
	\end{alignat}
\end{lem}

\begin{Proof}
Let $\tt \cq \tt_\cc$.
The upper bound in \cref{eq:ip:char:raw:all} implies \cref{eq:ip:char:pwr:short}:
if $\rr \ge \tfrac1{100} \nn$ and $\qq\lambda \ge 0$, then
\[
	\qq\lambda^\tt
\le
	\expb{ - \tfrac{201}{200} \tfrac{\rr}{\nn} \tt_\cc }
=
	\expb{ - \tfrac12 \rr \rbr{ \log \kk + \cc } \rbb{ 1 + \tfrac1{200} } }
=
	\rbr{ \ee^{\cc} \kk }^{-\rr(1+1/200)/2},
\]
using the inequality $1 - \xx \le \exp{ - \tfrac{201}{200} \xx }$, valid for all $\xx \ge \tfrac1{100}$.
Next, the lower bound in \cref{eq:ip:char:raw:all} implies that if $\rr \le \tfrac14 \nn$, then $\qq\lambda \ge 0$.
Thus, \cref{eq:ip:char:raw:all,eq:ip:char:raw:long}
implies \cref{eq:ip:char:pwr:long:upper}:
if $\rr \le \tfrac14 \nn$,~then%
\[
	0
\le
	\qq\lambda^\tt
\le
	\expb{ - \tfrac{2\rr}{\nn} \rbb{ 1 - \tfrac{\rr}{\nn} } \tt_\cc }
\le
	(\ee^\cc \kk)^{-\rr} \expb{ 2 \rr^2 \log \kk / \nn },
\]
using the inequality $1 - \xx \le \exp{ - \xx }$, valid for all $\xx \in \mbr$, and $\tt_\cc \le \nn \log \kk$, valid for $\kk$ large enough in terms of $\cc$.
For the corresponding lower bound \cref{eq:ip:char:pwr:long:lower}, we use the lower bound in \cref{eq:ip:char:raw:all} along with the inequality $1 - \xx \ge \exp{ - \xx - \xx^2 }$, valid for all $\xx \le \tfrac12$:
if $\rr \le \tfrac14 \kk$, then
\[
	\qq\lambda
\ge
	\expb{ 
	-	\tfrac{2\rr}{\nn} \rbb{ 1 + \tfrac{\rr}{\nn} }
	-	\tfrac{4\rr^2}{\nn^2} \rbb{ 1 + \tfrac{\rr}{\nn} }^2
	}
\ge
	\expb{ 
	-	\tfrac{2\rr}{\nn}
	-	\tfrac{7\rr^2}{\nn^2}
	}.
\]
Inequality~\cref{eq:ip:char:pwr:long:lower} follows by raising this to the power $\tt = \tt_\cc$ and using $\tt_\cc \le \nn \log \kk$.
\end{Proof}

We only use $\rr \asymp 1$ in the main term.
In particular, $\rr \asymp 1$ implies that $\rr \le \tfrac14 \nn$, so $\qq\lambda \ge 0$.

\begin{cor}
\label{res:ip:char:raw-dim:asymp}
	Suppose that $\lambda \vdash \nn$ is such that $\rr \cq \nn - \lambda_1 \asymp 1$.
	Then,
	\[
		\dd\lambda \qq\lambda^{\tt_\cc}
	=
		\tfrac1{\rr!} \dd*\lambda (\ee^{-\cc} \nn/\kk)^\rr
	\cdot
		\rbb{ 1 + \Oh{ \log \kk / \nn } }.
	\]
\end{cor}

\begin{Proof}
This is a simple consequence of \cref{res:prelim:rep:d_lambda} and \cref{eq:ip:char:pwr:long:upper,eq:ip:char:pwr:long:lower} from \cref{res:ip:char:pwr}.
\end{Proof}

\subsection{Control of the Error Term}
\label{sec:ip:error}

The first thing to check is that the error term in Teyssier's approximation lemma is small at the proposed mixing time $\tt_\cc$.
To do this, we need an upper bound on $\sqrt{ \mm_\lambda \dd\lambda } \qq\lambda^\tt$ at $\tt = \tt_\cc$;
see \cref{res:prelim:rep:mult-reweighted,res:ip:char:raw}.
We use these estimates to show that asymptotically all the mass is distributed over partitions $\lambda \vdash \nn$ with very long first row, namely $\nn - \lambda_1 \asymp 1$.

We remind the reader that throughout this section $\alpha \in (0, 1)$ is some fixed number,
	not depending on the underlying number $\nn$ of cards.
We assume that $\kk = \alpha \nn$ and ignore any non-integer issues.
We also let $\cc \in \mbr_+$ be arbitrary, but fixed, and take $\tt \cq \tt_\cc$.

\begin{prop}[Error Term]
\label{res:ip:error}
	Let $\alpha \in (0, 1)$ and $\kk = \alpha \nn$.
	Let $\cc \in \mbr$ and let $\tt \cq \tt_\cc$.
	Then,
	\[
		\LIMSUP{\toinf \MM} \,
		\LIMSUP{\toinf \nn} \,
		\sumt*{\lambda \vdash \nn : \lambda_1 < \nn-\MM}
		\sqrt{ \mm_\lambda \dd\lambda }
		\abs{\qq\lambda}^\tt
	=
		0.
	\]
\end{prop}

\begin{Proof}
	%
We start by using the bound
	on $\sqrt{ \mm_\lambda \dd\lambda }$ from \cref{res:prelim:rep:mult-reweighted}:
\[
	\sumt{\lambda \vdash \nn : \lambda_1 < \nn-\MM}
	\sqrt{ \mm_\lambda \dd\lambda } \abs{\qq\lambda}^\tt
\le
	\sumt{\rr=\MM+1}[\nn]
	\tfrac1{\rr!} (\nn \kk)^{\rr/2}
	\sumt{\lambda \vdash \nn : \lambda_1 = \nn-\rr}
	\dd*\lambda \abs{\qq\lambda}^\tt.
\]
Recall that $\nn\kk = \kk^2 / \alpha$.
Using Cauchy--Schwarz and $\sumt{\rho \in \ft \GG} \dd\rho^2 = \abs \GG$ for any $\GG$, we obtain
\[
	\sumt{\lambda \vdash \nn : \lambda_1 = \nn-\rr}
	\dd*\lambda
\le
	\sqrt{ \abs{ \bra{ \lambda \vdash \nn } } \rr! }
\le
	\sqrt{ 2^\rr \rr! }.
\]
Let
\(
	\qq\rr
\cq
	\max_{\lambda \vdash \nn : \lambda_1 = \nn-\rr}
	\qq\lambda
\vee
	0.
\)
We use this to handle $\lambda$ with $\qq\lambda \ge 0$:
\[
	\sumt{\lambda \vdash \nn : \lambda_1 < \nn-\MM, \: \qq\lambda \ge 0}
	\sqrt{ \mm_\lambda \dd\lambda } \abs{\qq\lambda}^\tt
\le
	\sumt{\rr=\MM+1}[\nn]
	(2/\alpha)^{\rr/2} \qq\rr^\tt \kk^\rr / \sqrt{\rr!}.
\]
We split this sum at $\rr = \tfrac1{100} \nn$ and apply \cref{eq:ip:char:pwr:short,eq:ip:char:pwr:long:upper} from \cref{res:ip:char:pwr} to the respective~cases.

First, consider $\rr \le \tfrac1{100} \nn$.
We apply \cref{eq:ip:char:pwr:long:upper} to bound $\qq\rr^\tt$ for $\lambda$ with $\qq\lambda \ge 0$:
\[
	\sumt{\rr=\MM+1}[\nn/100]
	(2/\alpha)^{\rr/2} \qq\rr^\tt \kk^\rr / \sqrt{\rr!}
\le
	\sumt{\rr=\MM+1}[\nn/100]
	\expb{ \rr \rbr{ \tfrac12 \log(2/\alpha) + \abs \cc } + \rr^2 \log \kk / \nn } / \sqrt{\rr!}.
\]
Let $\bb_\rr$ denote the $\rr$-th summand on the right-hand side.
Their ratios satisfy
\[
	\bb_{\rr+1} / \bb_\rr
&
=
	\expb{ \tfrac12 \log(2/\alpha) + \abs \cc + (2\rr+1) \log \kk / \nn } / \sqrt{\rr+1}
\\&
\le
	\expb{ \cc' + 5 \rr \log \nn / \nn } / \sqrt \rr
\eqqcolon
	\ff(\rr)
\Qwhere
	\cc'
\cq
	\tfrac12 \log(2/\alpha) + \abs \cc.
\]
The function $\ff$ is decreasing until $\tfrac52 \nn / \log \nn$ and increasing after this.
Thus,
\[
	\maxt{\MM \le \rr \le \nn/100}
	\ff(\rr)
\le
	\max\brb{ \ff(\MM), \: \ff(\tfrac1{100} \nn) }.
\]
Assume that $\MM \le \nn^{1/4}$.
Then, the following hold:
\begin{alignat*}{2}
	\ff(\MM)
&=
	\expb{ \cc' + 5 \MM \log \nn / \nn } / \sqrt \MM
\to
	0
&&\Qas
	\toinf \MM;
\\
	\ff(\tfrac1{100} \nn)
&=
	\expb{ \cc' + \tfrac14 \log \nn } / \sqrt{\nn/100}
\to
	0
&&\Qas
	\toinf \nn.
\end{alignat*}
Hence, the sum is dominated by a geometric sum with decay parameter $\tfrac12$.
Thus,
\[
&	\sumt{\rr=\MM+1}[\nn/100]
	\expb{ \rr \cc' + 2 \rr^2 \log \nn / \nn } / \sqrt{\rr!}
\le
	\expb{ \MM \cc' + 2 \MM^2 \log \nn / \nn } / \sqrt{\MM!}
\\&\qquad
\le
	\expb{ \MM \cc' + 2 \log \nn / \sqrt \nn } / \sqrt{\MM!}
\to
	0
\Qas
	\toinf \MM.
\]
We turn out attention to $\rr \ge \tfrac1{100} \nn$.
We apply \cref{eq:ip:char:pwr:short} to bound $\qq\rr^\tt$:
\[
&	\sumt{\rr=\nn/100}[\nn]
	(2/\alpha)^\rr \qq\rr^\tt \kk^\rr / \sqrt{\rr!}
\le
	\sumt{\rr=\nn/100}[\nn]
	\ee^{\rr \cc'}
	\kk^{\rr(1/2 - 1/400)}
	(\ee/\rr)^{\rr/2}
\\&\qquad
\le
	\sumt{\rr \ge \nn/100}
	\ee^{\rr \cc'}
	\kk^{\rr/2 - \rr/400}
	(100 \ee / \nn)^{\rr/2}
=
	\sumt{\rr \ge \nn/100}
	\ee^{\rr (\cc' + 3)}
	\kk^{-\rr/400}
=
	\oh1.
\]
Combining these two partial sums, we see that the sum over $\lambda$ with $\qq\lambda \ge 0$ tends to $0$.

It remains to consider $\lambda$ with $\qq\lambda \le 0$.
These terms can be bounded easily, using relatively crude bounds, to give an $\oh1$ contribution.
Analogous bounding is done by \textcite[Lemma~4.1]{T:limit-profile}---see their terms $S_1$, $S_2$ and $S_3$ there.
\end{Proof}

The above ideas are lifted almost verbatim from \cite{T:limit-profile}, with only minor adjustments.
In fact, the error term is small for all $\cc \in \mbr$ and $\tt_\cc(\nn, \kk) = \tt_{\cc'}(\nn, \nn)$ with $\cc' \cq \cc + \log \alpha \asymp 1$.
This, along with the fact that $\mm_\lambda \le \dd\lambda$, means that we can deduce \cref{res:ip:error} from \cite[Lemma~4.1]{T:limit-profile}.
We included the above proof for completeness, handling the most important cases (ie, small $\rr = \nn - \lambda_1$), deferring only the edge cases $\qq\lambda \le 0$ to \cite[Lemma~4.1]{T:limit-profile}.

We also felt it beneficial for the reader to have this as a `warm up', getting used to using these character ratios and dimension bounds before the main event---ie, \cref{res:ip:main}.

\subsection{Evaluation of the Main Term}
\label{sec:ip:main}

We analyse $\lambda$ with $\rr \asymp 1$ in the main term; the cost is an $\oh1$ error term, by \cref{res:ip:error}.
For such partitions, we can use the asymptotic bound 
from \cref{res:ip:char:raw-dim:asymp}.
Recall that
\[
	\Sigma \sim \Unif{\mfG = \symgr_\nn}
\Qand
	\Pi \sim \Unif{\mfK = \symgr_{\nn-\kk} \times \symgr_1 \times ... \times \symgr_1 \cong \symgr_{\nn-\kk}}
\quad
	\text{independently}.
\]

\begin{prop}[Main Term]
\label{res:ip:main}
	Let $\alpha \in (0, 1)$ and $\kk = \alpha \nn$.
	Let $\Sigma \sim \Unif\mfG$ and $\Pi \sim \Unif\mfK$.
	Let $\cc \in \mbr$ and $\tt \cq \tt_\cc$.
	Then,
	\[
		\LIMSUP{\toinf \MM}
		\LIMSUP{\toinf \nn} \,
		\absBB{
		\begin{multlined}[c]
			\exb{
				\absb{
					\exb{
						\sumt*{\lambda \vdash \nn : \lambda_1 \ge \nn - \MM}
						\dd\lambda \qq\lambda^\tt \chi_\lambda(\Sigma \Pi)
					\mid
						\Sigma
					}
				}
			}
		\qquad\\\qquad
		-	\DD*1[\Pois(\alpha+\ee^{-\cc})][\Pois(\alpha)]
		\end{multlined}
		}
	=
		0.
	\]
\end{prop}

\begin{Proof}[Skeleton of Proof]
We use the parametrisation $\beta = \ee^{-\cc} / \alpha$ in the results cited below.
In this parametrisation, the target of \cref{res:ip:main} is
$\DD1[\Pois\rbr{ \alpha(1+\beta) }][\Pois(\alpha)]$.

First, we abstractly formulate the main term in terms of numbers of fixed points of permutations.
This uses the estimate on $\dd\lambda \qq\lambda^{\tt_\cc}$ (\cref{res:ip:main:abst:char-est}) and evaluation of a representation-theoretic sum and polynomial (\cref{res:ip:main:abst:poly,res:ip:main:abst:pois-eval}), culminating in \cref{res:ip:main:abst:approx-combined}.

Next, we then split this number of fixed points into two parts (\cref{res:ip:main:abst:fixed-points}), resulting in \cref{res:ip:main:abst:approx-split}.
These numbers of fixed points can be approximated by Poisson random variables (\textit{\RV}s), for which simplifications and certain closed-form solutions exist.
Next, we evaluate the estimate of \cref{res:ip:main:abst:approx-split} under these approximations (\cref{res:ip:main:est:poisson-mgf,res:ip:main:est:given-L,res:ip:main:est:hypergeo-mgf}),
resulting in \cref{res:ip:main:est:main}.
This shows how the Poisson distance arises naturally.

It remains to evaluate the fixed-point laws precisely (\cref{res:ip:main:exact:fp-uncond,res:ip:main:exact:fp-cond}) and justify the approximations (\cref{res:ip:main:just}) to complete the proof.
\end{Proof}

\subsubsection{Abstract Formulation in Terms of Numbers of Fixed Points}

First, we use \cref{res:ip:char:raw-dim:asymp} to replace $\dd\lambda \qq\lambda^{\tt_\cc}$ by $\tfrac1{\rr!} \beta^\rr \dd*\lambda$ where $\beta = \ee^{-\cc} / \alpha$.

\begin{lem}
\label{res:ip:main:abst:char-est}
	Let $\alpha \in (0, 1)$ and $\kk = \alpha \nn$.
	Let $\Sigma \sim \Unif\mfG$ and $\Pi \sim \Unif\mfK$.
	Let $\cc \in \mbr$ and $\tt \cq \tt_\cc$.
	Write $\beta \cq \ee^{-\cc} / \alpha$.
	Then,
	\[
		\LIMSUP{\toinf \MM}
		\LIMSUP{\toinf \nn} \,
		\absBB{
		\begin{multlined}[c]
			\exb{
				\absb{
					\exb{
						\sumt*{\lambda \vdash \nn : \lambda_1 \ge \nn - \MM}
						\dd\lambda \qq\lambda^\tt \chi_\lambda(\Sigma \Pi)
					\mid
						\Sigma
					}
				}
			}
		\qquad\\\qquad
		-	\exb{ 
				\absb{ 
					\exb{ 
						\sumt{1 \le \rr \le \MM}
						\tfrac1{\rr!} \beta^\rr
						\sumt{\lambda \vdash \nn : \lambda_1 = \nn - \rr}
						\dd*\lambda \chi_\lambda(\Sigma \Pi)
					\mid
						\Sigma
					}
				}
			}
		\end{multlined}
		}
	=
		0.
	\]
\end{lem}

\begin{Proof}
The estimates of \S\ref{sec:ip:char:raw}, culminating in \cref{res:ip:char:raw-dim:asymp}, give
\[
	\dd\lambda \qq\lambda^\tt
=
	\tfrac1{\rr!} \dd*\lambda \beta^\rr
\cdot
	\rbb{ 1 + \Oh{ \tfrac1\kk \log \nn } }
\Qwhen
	\rr
\cq
	\nn - \lambda_1
\asymp
	1.
\]
This bound is analogous to that obtained by \textcite[Equation~(4.2)]{T:limit-profile}, as $\kk \asymp \nn$.
The same simple, albeit somewhat technical, argument as they give proves our claim.
In essence, it is just controlling an error sum using some fairly crude bounds.
\end{Proof}

We could use linearity of the expectation to pass the inner $\ex\cdot$ through the sums, giving $\ex{ \chi_\lambda(\Sigma \Pi) \mid \Sigma }$.
However, evaluating this goes back to the spherical functions-type approach.
Rather, we evaluate the inner sums,
	first over $\lambda$ with fixed first row $\lambda_1 = \nn-\rr$,
	then over $\rr$.
Only after having done these computations do we take the expectation.


To this end, we introduce a collection $\bra{ \TT_\rr }_{\rr \in \mbn}$ of polynomials:
\[
	\TT_\rr(\zz)
\cq
	\sumt{\ii=0}[\rr]
	\binomt{\zz}{\rr-\ii}
	(-1)^\ii / \ii!
\Qfor
	\rr, \zz \in \mbn.
\]
\textcite[Lemma~4.3]{T:limit-profile} proves the following representation-theoretic result. It follows from the Murnaghan--Nakayama rule for calculating the character along with a clever, non-standard choice of basis.
Let $\symgr_{\nn, \rr}$ denote the set of permutations $\sigma \in \symgr_\nn$ with all cycles of length at most $\rr$.
For a permutation $\sigma$, let $\Fix \sigma$ denote its number of fixed points:
\[
	\Fix \sigma
\cq
	\abs{ \bra{ \ii \mid \sigma(\ii) = \ii } }.
\]

\begin{lem}[{\textcite[Lemma~4.3]{T:limit-profile}}]
\label{res:ip:main:abst:poly}
	Let $\rr \in \mbn$ and let $\sigma \in \symgr_\nn$ be a permutation with some cycle of length greater than $\rr$, ie $\sigma \in \symgr_\nn \setminus \symgr_{\nn, \rr}$.
	Then,
	\[
		\tfrac1{\rr!}
		\sumt{\lambda \vdash \nn : \lambda_1 = \nn-\rr}
		\dd*\lambda \chi_\lambda(\sigma)
	=
		\TT_\rr\rbr{\Fix \sigma}.
	\]
\end{lem}

The above lemma only applies for permutations $\sigma$ with minimal cycle-length at least $\rr$; if the minimal length is at least $\MM$, then this can be used for all $\rr \le \MM$ simultaneously.
Later, we need to average over $\sigma \in \symgr_\nn$; the minimal length is larger than $\MM$ \whp (\asinf \nn), so we are able to exclude cases for which this lemma does not apply.

\begin{lem}
\label{res:ip:main:abst:short-cycles}
	Let $\MM \in \mbn$ and $\beta \in (0, \infty)$.
	Then,
	\[
		\LIMSUP{\toinf \nn} \,
		\abs {\symgr_\nn}^{-1}
		\sumt{\sigma \in \symgr_{\nn, \MM}}
		\sumt*{\lambda \vdash \nn : \lambda_1 \ge \nn - \MM}
		\rbb{
			\dd*\lambda \abs{ \chi_\lambda(\sigma) }
		+	\abs{ \beta^\rr \TT_\rr(\Fix \sigma) }
		}
	=
		0.
	\]	
\end{lem}

\begin{Proof}
This follows from some rough bounds, as shown by \textcite[Proof of Lemma 4.2]{T:limit-profile}---there, the sum in question is called $S_{\symgr_{\nn, \ell}}$, with their $\ell$ replacing our $\MM$.
\end{Proof}

This lemma allows us to, in essence, apply \cref{res:ip:main:abst:poly} for all $\sigma \in \symgr_\nn$, rather than just those in $\symgr_\nn \setminus \symgr_{\nn, \MM}$:
	the impact of the $\sigma \in \symgr_{\nn, \MM}$ is insignificant on the overall sum.

We must sum over these polynomials, each weighted by $\beta^\rr$.
We think of $M$ as ``arbitrarily large, but fixed''.
Next, we approximate $\sumt{1 \le \rr \le \MM}$ by $\sumt{\rr \ge 1}$.
This is done by \textcite[Lemma~4.4]{T:limit-profile} using elementary, if somewhat technical, arguments.

\begin{lem}[{\cite[Lemma~4.4]{T:limit-profile}}]
\label{res:ip:main:abst:r-inf}
	Let $\NN \in \mbn$ and $\beta \in (0, \infty)$.
	Then,
	\[
		\LIMSUP{\toinf \MM} \,
		\absb{
			\sumt{1 \le \rr \le \MM}
			\beta^\rr
			\TT_\rr(\NN)
		-	\sumt{\rr \ge 1}
			\beta^\rr
			\TT_\rr(\NN)
		}
	=
		0.
	\]
\end{lem}

Next, we must evaluate the sum $\sumt{\rr\ge1}$.
This is done by \textcite[Proposition~4.5]{T:limit-profile} using a simple change of variables and swapping of summation order.

\begin{lem}[{\cite[Proposition~4.5]{T:limit-profile}}]
\label{res:ip:main:abst:pois-eval}
	Let $\NN \in \mbn$ and $\beta \in (0, \infty)$.
	Then,
	\[
		\sumt{\rr \ge 1}
		\tfrac1{\rr!} \beta^\rr
		\TT_\rr(\NN)
	=
		\ff_\beta(\NN)
	\Qwhere
		\ff_\beta(\NN)
	\cq
		\ee^{-\beta} (1 + \beta)^\NN - 1.
	\]
\end{lem}

Combining \cref{res:ip:main:abst:char-est,res:ip:main:abst:poly,res:ip:main:abst:short-cycles,res:ip:main:abst:r-inf,res:ip:main:abst:pois-eval} gives the following approximation to the $\ell_1$ distance.

\begin{cor}
\label{res:ip:main:abst:approx-combined}
	Let $\alpha \in (0, 1)$ and $\kk = \alpha \nn$.
	Let $\Sigma \sim \Unif\mfG$ and $\Pi \sim \Unif\mfK$.
	Let $\cc \in \mbr$ and $\tt \cq \tt_\cc$.
	Write $\beta \cq \ee^{-\cc} / \alpha$.
	Then,
	\[
		\LIMSUP{\toinf \MM}
		\LIMSUP{\toinf \nn} \,
		\absBB{
		\begin{multlined}[c]
			\exb{
				\absb{
					\exb{
						\sumt*{\lambda \vdash \nn : \lambda_1 \ge \nn - \MM}
						\dd\lambda \qq\lambda^\tt \chi_\lambda(\Sigma \Pi)
					\mid
						\Sigma
					}
				}
			}
		\qquad\\\qquad
		-	\exb{
				\absb{
					\ee^{-\beta}
					\exb{ (1 + \beta)^{\Fix(\Sigma\Pi)} \mid \Sigma }
				-	1
				}
			}
		\end{multlined}
		}
	=
		0.
	\]
\end{cor}

Recall that $\mfK = \symgr_\mm \times \symgr_1 \times ... \times \symgr_1 \subseteq \symgr_\nn = \mfG$, so the averaging over $\mfK$ corresponds to uniformising the $\mm$ cards.
If $\mfK = \bra{\id}$ is the trivial group---ie, $\kk = \nn$, so all cards are followed---then the averaging over $\Pi \sim \Unif{\mfK}$ does nothing and the evaluation becomes relatively simple.
A quantitative version of the well-known fact that $\Fix \Sigma \approx \Pois(1)$ gives a fairly simple proof; see \cite[Lemma~4.6]{T:limit-profile}.
It is much more challenging in general.

The fundamental idea is to break up $\Fix(\Sigma \Pi)$ into two parts:
	the fixed points amongst the first $\mm \cq \nn-\kk$ indices and the remainder.
We then estimate various expectations via approximating \RVs by Poisson \RVs and using concentration of a certain hypergeometric \RV.

Below, all indices are assumed implicitly to be in $[\nn] = \bra{1, ..., \nn}$.

\begin{lem}
\label{res:ip:main:abst:fixed-points}
	Let $\sigma \in \symgr_\nn$ and
	\(
		\LL_\sigma
	\cq
		\abs{ \bra{ \ii \le \mm \mid \sigma^{-1}(\ii) \le \mm } }.
	\)
	Let $\Pi \sim \Unif{\mfK}$.
	Then,
	\[
		\Fix(\sigma\Pi)
	=^d
		\abs{ \bra{ \ii \le \LL_\sigma \mid \Pi(\ii) = \ii } }
	+	\abs{ \bra{ \ii > \mm \mid \sigma(\ii) = \ii } }.
	\]
\end{lem}

\begin{Proof}
Let $\mcl_\sigma \cq \bra{ \ii \le \mm \mid \sigma^{-1}(\ii) \le \mm}$.
Any $\pi \in \mfK$ stabilises indices $\ii > \mm$.
Thus,
\[
	\Fix(\sigma \pi)
&
=
	\abs{
		\bra{ \ii \le \mm \mid \pi(\ii) = \sigma^{-1}(\ii) }
	\mathrel{\dot \cup}
		\bra{ \ii >   \mm \mid \sigma(\ii) = \ii }
	}
\\&
=
	\abs{ \bra{ \ii \le \mm \mid \pi(\ii) = \sigma^{-1}(\ii) } \cap \mcl_\sigma }
+	\abs{ \bra{ \ii > \mm \mid \sigma(\ii) = \ii } }.
\]
Now, $\sigma$ is a fixed permutation.
Thus, by symmetry, the law of $\Fix(\sigma \Pi)$ is unaffected by the particular choice of either $\mcl_\sigma$ or $\sigma^{-1}$.
We may thus assume that $\mcl_\sigma = [\LL_\sigma]$ and $\sigma^{-1}(\ii) = \ii$ for all $\ii \le \LL_\sigma$; this does not affect the law of $\Fix(\sigma \Pi)$.
This proves the lemma.
\end{Proof}

This lemma deconstructs $\Fix(\Sigma \Pi)$ into two conditionally independent parts given $\LL_\Sigma$.

\begin{cor}
\label{res:ip:main:abst:approx-split}
	Let $\alpha \in (0, 1)$ and $\kk = \alpha \nn$.
	Let $\Sigma \sim \Unif{\symgr_\nn}$ and $\Pi \sim \Unif{\symgr_\mm \times \symgr_1 \times \cdots \times \symgr_1 \subseteq \symgr_\nn}$.
	Let $\cc \in \mbr$ and $\tt \cq \tt_\cc$.
	Write $\beta \cq \ee^{-\cc} / \alpha$.
	Then,
	\begin{gather*}
		\DD1[\tt]
	=
		\exb{ \absb{ \ee^{-\beta} \exb{ (1 + \beta)^{\XX + \YY} \mid \LL_\Sigma } - 1 } }
	\quad
		\text{where}
	\\
		\XX
	\cq
		\abs{ \bra{ \ii \le \LL_\Sigma \mid \Pi(\ii) = \ii } }
	\Qand
		\YY
	\cq
		\abs{ \bra{ \ii > \mm \mid \Sigma(\ii) = \ii } }.
	\end{gather*}
	Moreover, $\XX$ and $\YY$ are conditionally independent given $\LL_\Sigma = \abs{ \bra{ \ii \le \mm \mid \Sigma^{-1}(\ii) \le \mm } }$.
\end{cor}

We can, and do, calculate the law of each of $\XX$ and $\YY$ explicitly, given the value of $\LL_\Sigma$. These are not amenable to closed-form simplification, however.
The laws we calculate are each approximately Poisson, given $\LL_\Sigma$, with an explicit parameter:
\[
	\XX
\approx
	\Pois\rbb{\tfrac \LL\nn}
\Qand
	\YY
\approx
	\Pois\rbb{1 - \tfrac{\mm-\LL}{\nn-\mm}}
\Quad{given}
	\LL_\Sigma = \LL.
\]
Further, the law of $\LL_\Sigma$ is \emph{exactly} hypergeometric.
	A hypergeometric distribution $\HG(\AA, \BB, \nn)$ is one that ``describes the probability of a given number of successes in $\AA$ draws, without replacement, from a finite population of size $\nn$ that contains exactly $\BB$ objects with that feature''.
	$\LL_\Sigma$ fits this precisely with $\AA = \BB = \mm$:
		$\LL_\Sigma \sim \HG(\mm, \mm, \nn)$.
We record this now.

\begin{lem}
\label{res:ip:main:abst:hypergeo}
	Let $\Sigma \sim \Unif{\symgr_\nn}$.
	Then,
	\[
		\LL_\Sigma
	=
		\abs{ \bra{ \ii \le \mm \mid \Sigma^{-1}(\ii) \le \mm } }
	\sim
		\HG(\mm, \mm, \nn).
	\]
\end{lem}

Using the Poisson approximations and known concentration results for hypergeometrics%
	---they concentrate more strongly than the corresponding ``with replacement'' Binomial distribution $\Bin(\AA, \BB/\nn)$---%
we are able to well-approximate the relative $\ell_1$ distance $\DD1[\tt]$.
We do this first, then determine precisely the laws of
\[
	\abs{ \bra{ \ii \le \LL_\sigma \mid \Pi(\ii) = \ii } }
\Qand
	\abs{ \bra{ \ii > \mm \mid \sigma(\ii) = \ii } },
\]
showing that they are approximately Poisson.
Finally, we briefly justify the replacement of these by their approximations in the $\ell_1$-distance formula.

\subsubsection{Estimation Using Poisson Approximations}

We first estimate the $\ell_1$ distance from \cref{res:ip:main:abst:approx-split} when replacing $\XX$ and $\YY$ with their conditionally-independent Poisson approximations.

\begin{prop}
\label{res:ip:main:est:main}
	Let $\alpha \in (0, 1)$ and $\kk = \alpha \nn$.
	Let $\LL \sim \HG(\mm, \mm, \nn)$. Let $\XX_0 \sim \Pois(\tfrac \LL\mm)$ and $\YY_0 \sim \Pois(1 - \tfrac{\mm-\LL}{\nn-\mm})$ independently given $\LL$.
	Let $\cc \in \mbr$ and $\tt \cq \tt_\cc$. Write $\beta \cq \ee^{-\cc} / \alpha$.
	Then,
	\[
		\LIMSUP{\toinf \nn} \,
		\absb{
			\exb{ \absb{ \ee^{-\beta} \exb{ (1 + \beta)^{\XX_0 + \YY_0} \mid \LL } - 1 } }
		-	\DD*1[\Pois\rbr{ \alpha (1 + \beta) }][\Pois(\alpha)]
		}
	=
		0.
	\]
\end{prop}

We prove this via a sequence of lemmas.
%
The first corresponds primarily to the inner expectation. It is straightforward since there are no absolute value signs to cause difficulties.
We omit its simple proof: it goes via the usual moment-generating function (\textit{\mgf}) of the~Poisson.

\begin{lem}
\label{res:ip:main:est:poisson-mgf}
	Let $\zz \in (0, \infty)$ and $\ZZ_0 \sim \Pois(\zz)$.
	Let $\gamma \in (0, \infty)$.
	Then,
	\[
		\exb{ (1 + \gamma)^{\ZZ_0} }
	=
		\exp{ \gamma \zz }.
	\]
\end{lem}

The next lemma is far more challenging technically.
We do not evaluate it exactly, but rather only up to a term which will be small when we average over $\LL \sim \HG(\mm, \mm, \nn)$, the mean of which is $\mm^2/\nn$.
This is due to the fact that the averaging is being done \emph{outside} the absolute value, leading to \TV-type expressions and other technical hurdles.

\begin{lem}
\label{res:ip:main:est:given-L}
	Let $\LL \in \mbn$ with $0 \le \mm - \LL \le \nn - \mm$.
	Let $\lambda_0 = 1 - \tfrac{\mm-\LL}{\nn-\mm}$ and $\YY_0 \sim \Pois(\lambda_0)$; let $\lambda_\star \cq 1 - \tfrac \mm\nn = \tfrac \kk\nn$, ie the value of $\lambda_0$ when $\LL = \mm^2/\nn$.
	Let $\beta \in (0, \infty)$.
	Then,
	\[
	&	\absb{
			\exb{ \absb{ \ee^{ - \beta + \beta \LL/\mm } (1 + \beta)^{\YY_0} - 1 } }
		-	\DD*1[\Pois\rbr{\lambda_\star(1+\beta)}][\Pois(\lambda_\star)]
		}
	\\&\qquad
	\le
		\ee^{\beta^2}
		\rbb{
			\expb{ \beta \tfrac \nn{\nn-\mm} \absb{ \tfrac \LL\mm - \tfrac \mm\nn } }
		-	1
		}^2
	+	\ee^{\beta^2}
		\expb{ \beta \tfrac \nn{\nn-\mm} \absb{ \tfrac \LL\mm - \tfrac \mm\nn } }
		\rbb{
			\expb{ \beta^2 \abs{ \tfrac \LL\mm - \tfrac \mm\nn } }
		-	1
		}.
	\]
\end{lem}

\begin{Proof}
	%
Let $\YY_\star \sim \Pois(\lambda_\star)$.
We `pivot' $\LL$ around $\mm^2/\nn$ to convert $\YY_0$ into $\YY_\star$:
\[
&	\exb{ \absb{ \ee^{ - \beta + \beta \LL/\mm } (1 + \beta)^\YY - 1 } }
-	\exb{ \absb{ \ee^{-\beta\lambda_\star} (1 + \beta)^{\YY_\star} - 1 } }
\\&\qquad
=
	\ex{
		\absb{
			\ee^{-\beta\lambda_\star}
			(1 + \beta)^{\YY_*}
		\cdot
			\ee^{ \beta \rbr{ \LL/\mm - \mm/\nn } }
			(1 + \beta)^{\YY_0 - \YY_\star}
		-	1
		}
	-	\absb{
			\ee^{-\beta\lambda_\star} (1 + \beta)^{\YY_\star} - 1
		}
	}
\\&\qquad
\le
	\exb{
		\ee^{-\beta\lambda_\star}
		(1 + \beta)^{\YY_\star}
		\absb{ \ee^{\beta(\LL/\mm - \mm/\nn)} (1 + \beta)^{\YY_0 - \YY_\star} - 1 }
	}.
\]
The expectation about which we are pivoting is easy to calculate explicitly:
\[
	\exb{ \absb{ \ee^{-\beta\lambda_\star} (1 + \beta)^{\YY_\star} - 1 } }
&
=
	\sumt{\yy\ge0}
	\tfrac1{\yy!} \ee^{-\lambda_\star} \lambda_\star^\yy
	\absb{ \ee^{-\beta \lambda_\star} (1 + \beta)^\yy - 1 }
\\&
=
	\sumt{\yy\ge0}
	\absb{ \tfrac1{\yy!} \ee^{-\lambda_\star(1+\beta)} \rbb{ \lambda_\star(1+\beta) }^\yy - \tfrac1{\yy!} \ee^{-\lambda_\star} \lambda_\star^\yy }
\\&
=
	\DD*1[\Pois\rbr{\lambda_\star(1+\beta)}][\Pois(\lambda_\star)].
\]

It remains to bound the remaining `error'.
Let $\LL' \cq \LL - \mm^2/\nn$.
First, note that
\[
	\lambda_0 - \lambda_\star
=
	\frac{\mm}{\nn} - \frac{\mm-\LL}{\nn-\mm}
=
	\frac{ \mm(\nn - \mm) - \nn(\mm - \LL) }{\nn(\nn-\mm)}
=
	\frac{ \LL - \mm^2 / \nn }{\nn - \mm}
=
	\frac{\LL'}{\nn - \mm}.
\]
We now use a simple coupling of two Poisson \RVs:
	if $0 < \aa < \bb$, then $\Pois(\aa) + \Pois(\bb-\aa) \sim \Pois(\bb)$,
	if $\Pois(\aa)$ and $\Pois(\bb-\aa)$ are independent.
Thus,
\(
	\abs{ \YY_0 - \YY_\star }
\sim
	\Pois(\abs{\lambda_0 - \lambda_\star}).
\)
This is not independent of $\YY_\star$, however, so we apply Cauchy--Schwarz to separate the parts:
\[
&	\exb{
		\ee^{-\beta\lambda_\star}
		(1 + \beta)^{\YY_\star}
		\absb{ \ee^{\beta\LL'} (1 + \beta)^{\YY_0 - \YY_\star} - 1 }
	}{}^2
\\&\qquad
\le
	\exb{ \ee^{-2\beta\lambda_\star} (1 + \beta)^{2\YY_\star} }
\cdot
	\exb{ \rbb{ \ee^{\beta\abs{\LL'}/\mm} (1 + \beta)^\ZZ - 1 }{}^2 }
\]
where $\ZZ \sim \Pois\rbr{ \abs{\lambda_0 - \lambda_\star} } = \Pois\rbr{ \abs{\LL'} (\nn-\mm) }$.
The first term is straightforward to calculate and bound using \cref{res:ip:main:est:poisson-mgf} and the fact that $\YY_\star \sim \Pois(\lambda_\star)$ and $\lambda_\star \le 1$:
\[
	\exb{ \ee^{-2\beta\lambda_\star} (1 + \beta)^{2\YY_\star} }
=
	\ee^{\beta^2\lambda_\star}
	\exb{ \ee^{-(2\beta+\beta^2)\lambda_\star}(1 + 2\beta + \beta^2)^{\YY_\star} }
=
	\ee^{\beta^2 \lambda_\star}
\le
	\ee^{\beta^2}.
\]
Next, we expand the square in the second term and use \cref{res:ip:main:est:poisson-mgf} again:
\[
&	\exb{ \rbb{ \ee^{\beta\abs{\LL'}/\mm} (1 + \beta)^\ZZ - 1 }^2 }
=
	\ee^{2\beta\abs{\LL'}/\mm}
	\exb{ (1 + \beta)^{2\ZZ} }
-	2 \ee^{\beta\abs{\LL'}/\mm}
	\exb{ (1 + \beta)^{\ZZ} }
+	1
\\&\qquad
=
	\ee^{2\beta\abs{\LL'}/\mm} \ee^{(2\beta+\beta^2)\abs{\lambda_0 - \lambda_\star}}
-	2 \ee^{\beta\abs{\LL'}/\mm} \ee^{\beta\abs{\lambda_0 - \lambda_\star}}
+	1.
\\
\intertext{Plugging in the expressions for $\lambda_0$ and $\lambda_\star$, we find that this is equal to}
&\qquad
\mathrel{\phantom{=}}
	\ee^{2\beta\abs{\LL'}(1/\mm + 1/(\nn-\mm))} \ee^{\beta^2 \abs{\LL'}/(\nn-\mm)}
-	2 \ee^{\beta\abs{\LL'}(1/\mm + 1/(\nn-\mm))}
+ 	1
\\&\qquad
=
	\rbb{ \ee^{\beta \abs{\LL'} (1/\mm + 1/(\nn-\mm))} - 1 }^2
+	\ee^{2\beta\abs{\LL'}(1/\mm+1/(\nn-\mm))}
	\rbb{ \ee^{\beta^2\abs{\LL'}/\mm} - 1 }.
\]

The claim now follows immediately from all these estimates and calculations.
\end{Proof}

The `error bound' in the previous lemma is, admittedly, pretty confusing.
In particular, it is not `small' unless $\LL$ is `close to' $\mm^2/\nn$.
We do not try to explicitly using the distribution $\LL \sim \HG(\mm, \mm, \nn)$.
Rather, we use the fact that $\LL$ concentrates very well around its mean $\mm^2/\nn$; further, $\LL/\mm \in [0, 1]$, so is uniformly bounded. This will allow us to replace $\LL/\mm$ by $\mm/\nn + \oh{\mm/\nn}$ at only an $\oh1$ error \emph{even in the \mgfs}.

Hypergeometric distributions,
	which correspond to sampling without replacement,
concentrate more strongly than their Binomial counterparts,
	which correspond to sampling with replacement.
This can be made formal by a simple coupling.
The well-known Binomial concentration bound is sufficient for our application.

\begin{lem}
\label{res:ip:main:est:hypergeo-mgf}
	Let $\LL \sim \HG(\mm, \mm, \nn)$.
	Let $\gamma, \omega \in (0, \infty)$ with $\gamma \omega \le \mm$.
	Then,
	\[
		\absb{ \exb{ \expb{ \gamma \abs{ \LL/\mm - \mm/\nn } } } - 1 }
	\le
		2 \gamma \omega / \mm + 2 \exp{ - 2 \omega^2 / \mm }.
	\]
\end{lem}

\begin{Proof}
The claimed bound follows easily from the standard hypergeometric tail bound
\[
	\pr{ \abs{ \LL/\mm - \mm/\nn } \ge \omega/\mm }
\le
	2 \exp{ - 2 \omega^2 \mm },
\]
the uniformity $\LL/\mm \in [0, 1]$ and the inequality $\ee^x - 1 \le 2x$ for $x \in [0, 1]$.
Indeed,
\[
	\absb{ \exb{ \expb{ \gamma \abs{ \LL/\mm - \mm/\nn } } } - 1 }
&
\le
	\rbb{ \exp{ \gamma \omega / \mm } - 1 }
+	\pr{ \abs{ \LL/\mm - \mm/\nn } \ge \omega/\mm }
\\&
\le
	2 \gamma \omega / \mm + 2 \exp{ - 2 \omega^2 / \mm }.
\qedhere
\]
\end{Proof}

These results combine easily to estimate the distance under the Poisson approximation.

\begin{Proof}[Proof of \cref{res:ip:main:est:main}]
The proposition follows from applying \cref{res:ip:main:est:hypergeo-mgf} to each of the terms in \cref{res:ip:main:est:given-L},
noting that
\(
	\abs{ \tfrac \LL\mm - \tfrac \mm\nn } \le 1
\Qand
	\lambda_\star = 1 - \tfrac \mm\nn = \tfrac \kk\nn = \alpha.
\)
\end{Proof}

\subsubsection{Exact Calculation of Fix-Point Laws}

We need to determine the law of $\LL_\Sigma$ and understand the number of fixed points of a permutation amongst a certain collection of indices.
We analyse the number of fixed point~first.

\begin{lem}
\label{res:ip:main:exact:fp-uncond}
	Let $\LL \in \mbn$ with $\LL \le \mm$ and let $\Pi \sim \Unif{\symgr_\mm}$.
	Then,
	\[
		\pr{ \abs{ \bra{ \ii \le \LL \mid \Pi(\ii) = \ii } } = \rr }
	=
		\frac1{\rr!}
		\sumd{0 \le \ell \le \LL-\rr}
		\frac{(-1)^\ell}{\ell!}
		\frac{(\LL)_{\rr+\ell}}{(\mm)_{\rr+\ell}}.
	\]
	In particular, if $\LL \gg 1$, then
	\(
		\abs{ \bra{ \ii \le \LL \mid \Pi(\ii) = \ii } }
	\approx
		\Pois(\frac{\LL}{\mm}).
	\)
\end{lem}

\begin{Proof}
This result is well known when $\LL = \mm$.
General $\LL$ requires some adjustment.
Write
\[
	\Fix_\LL \pi \cq \abs{ \bra{ \ii \le \LL \mid \pi(\ii) = \ii } }
\Qand
	\FF_\ii \cq \bra{ \pi \in \symgr_\mm \mid \pi(\ii) = \ii }
\]
for
	the number of fixed points of $\pi \in \symgr_\mm$ amongst the first $\LL$ indices
and
	the set of $\pi \in \symgr_\mm$ with $\ii \in [\mm]$ as a fixed point,
respectively.
Then,
\[
	\bra{ \pi \in \symgr_\mm \mid \Fix_\LL \pi = 0 }
=
	\cap_{\ii \le \LL} \FF_\ii^c
=
	\rbr{ \cup_{\ii \le \LL} \FF_\ii }^c.
\]
By the inclusion--exclusion principle,
\[
	\abs{ \cup_{\ii \le \LL} \FF_\ii }
=
	\sumt{\emptyset \ne \JJ \subseteq [\LL]}
	(-1)^{\abs \JJ + 1}
	\abs{ \cap_{\jj \in \JJ} \FF_\jj }.
\]
By symmetry,
\(
	\abs{ \cap_{\jj \in \JJ} \FF_\jj }
=
	(\mm - \abs \JJ)!.
\)
Hence,
\[
	\abs{ \cup_{\ii \le \LL} \FF_\ii }
&
=
	\sumt{1 \le \ell \le \LL}
	\binomt \LL\ell
	(-1)^{\ell+1}
	(\mm - \ell)!
=
	\mm!
	\sumt{1 \le \ell \le \LL}
	\tfrac1{\ell!} (-1)^{\ell+1} (\LL)_\ell/(\mm)_\ell.
\]
In particular, combining these gives
\[
	\abs{ \bra{ \pi \in \symgr_\mm \mid \Fix_\LL \pi = 0 } }
=
	\mm!
	\sumd{0 \le \ell \le \LL}
	\frac{(-1)^\ell}{\ell!}
	\frac{(\LL)_\ell}{(\mm)_\ell}.
\]
This recovers the well-known approximately-$\Pois(1)$ formula when $\LL = \mm$.

We now generalise this to
$\bra{ \pi \in \symgr_\mm \mid \Fix_\LL \pi = \rr }$
for $\rr \ge 1$.
We have
\[
	\bra{ \pi \in \symgr_\mm \mid \Fix_\LL \pi = \rr }
=
	\cup_{\II \subseteq [\LL] : \abs \II = \rr}
	\rbb{
		\rbr{ \cap_{\ii \in \II} \FF_\ii }
	\cap
		\rbr{ \cap_{\ii \notin \II} \FF_\ii^c }
	}
\]
The event $\cap_{\ii \in \II} \FF_\ii$ forces the indices in $\II \subseteq [\LL]$ to be fixed points.
The remaining indices of $[\LL]$, ie $[\LL] \setminus \II$, must be non-fixed points, whilst those in $[\mm] \setminus [\LL]$ are unrestricted.
There are $\binom \mm\rr$ choices for $\II$; all are symmetric, so we may assume that $\II = \bra{\LL-\rr+1, ..., \LL}$.
Thus,
\[
&	\abs{ \bra{ \pi \in \symgr_\mm \mid \Fix_\LL \pi = \rr } }
=
	\binom \LL\rr
	(\mm-\rr)!
	\sumd{0 \le \ell \le \LL-\rr}
	\frac{(-1)^\ell}{\ell!}
	\frac{(\LL-\rr)_\ell}{(\mm-\rr)_\ell}.
\\&\qquad
=
	\frac{(\LL)_\rr}{\rr!}
	\frac{\mm!}{(\mm)_\rr}
	\sumd{0 \le \ell \le \LL-\rr}
	\frac{(-1)^\ell}{\ell!}
	\frac{(\LL-\rr)_\ell}{(\mm-\rr)_\ell}
=
	\frac{\mm!}{\rr!}
	\sumd{0 \le \ell \le \LL-\rr}
	\frac{(-1)^\ell}{\ell!}
	\frac{(\LL)_{\rr+\ell}}{(\mm)_{\rr+\ell}}.
\]
The first part lemma now follows since $\Pi \sim \Unif{\symgr_\mm}$, so
\[
	\pr{ \Fix_\LL \pi = \rr }
=
	\abs{ \bra{ \pi \in \symgr_\nn \mid \Fix_\LL \pi = \rr } } / \mm!.
\]

This shows that the distribution of $\Fix_\LL \Pi$ is approximately $\Pois(\LL/\mm)$ when $1 \ll \LL \le \mm$.
This is well known when $\LL = \mm$ and follows more generally because $(\LL-\rr)_\ell / (\mm-\rr)_\ell \approx (\LL/\mm)^\ell$ when $\rr, \ell \asymp 1$.
Asymptotically all the mass of the sum comes from $\ell \asymp 1$ and $\Fix_\LL \Pi \le \Fix(\Pi) \approx \Pois(1)$, so the important case is $\rr, \ell \asymp 1$.
\end{Proof}

Exactly the same argument can be translated into the event for $\Sigma$.
However, there is a slight subtlety:
	the \RV $\LL_\Sigma$ is correlated with $\abs{ \bra{ \ii > \mm \mid \Sigma(\ii) = \ii } }$.
Knowing the value of $\LL_\Sigma$ determines the number of indices $\ii > \mm$ which are mapped to $\Sigma(\ii) \le \mm$, which cannot be fixed points.
Given this, it follows an analogous distribution to that of the previous lemma.

\begin{lem}
\label{res:ip:main:exact:fp-cond}
	Let $\LL \in \mbn$ with $0 \le \mm - \LL \le \nn - \mm$ and $\Sigma \sim \Unif{\symgr_\nn}$.
	Then,
	\[
	&	\pr{ \abs{ \bra{ \ii > \mm \mid \Sigma(\ii) = \ii } } = \rr \mid \LL_\Sigma = \LL }
	=
		\frac1{\rr!}
		\sumd{0 \le \ell \le \nn-2\mm+\LL - \rr}
		\frac{(-1)^\ell}{\ell!}
		\frac{(\nn-2\mm+\LL)_{\rr+\ell}}{(\nn-\mm)_{\rr+\ell}}.
	\]
	In particular, if $\nn-2\mm+\LL \gg 1$, then
	\(
		\abs{ \bra{ \ii > \mm \mid \Sigma(\ii) = \ii } }
	\approx
		\Pois(1 - \frac{\mm-\LL}{\nn-\mm})
	\)
	given $\LL_\Sigma = \LL$.
\end{lem}

\begin{Proof}
This follows analogously to the above after taking care of the conditioning.
We have
\[
	\abs{ \bra{ \ii \le \mm \mid \Sigma^{-1}(\ii) \le \mm } }
=
	\LL
&
\iff
	\abs{ \bra{ \ii > \mm \mid \Sigma(\ii) \le \mm } }
=
	\mm - \LL
\\&
\iff
	\abs{ \bra{ \ii > \mm \mid \Sigma(\ii) > \mm } }
=
	(\nn - \mm) - (\mm - \LL).
\]
We are choosing uniformly without replacement $\nn-\mm$ indices $\ii > \mm$ and conditioning that $(\nn - \mm) - (\mm - \LL) = \nn - 2\mm + \LL$ are mapped to indices $\jj > \mm$%
	---which is a necessary condition for $\ii$ to be a fixed point, ie $\ii = \jj$.
The number of fixed points amongst indices $\ii > \mm$ is then approximately $\Pois(\frac{\nn-2\mm+\LL}{\nn-\mm}) = \Pois(1 - \frac{\mm-\LL}{\nn-\mm})$ in the precise sense of \cref{res:ip:main:exact:fp-uncond}.
\end{Proof}

The above proof also shows that the two parts are independent given $\LL_\Sigma$.
Recall that
\[
	\XX
\cq
	\abs{ \bra{ \ii \le \LL_\Sigma \mid \Pi(\ii) = \ii } }
\Qand
	\YY
\cq
	\abs{ \bra{ \ii > \mm \mid \Sigma(\ii) = \ii } }.
\]

\begin{cor}
\label{res:ip:main:exact:indep}
	The \RVs
	\(
		\XX
	\)
	and
	\(
		\YY
	\)
	are conditionally independent given $\LL_\Sigma$.
	Moreover,
	\[
		\XX
	\approx
		\Pois\rbb{\tfrac \LL\mm}
	\Qand
		\YY
	\approx
		\Pois\rbb{1 - \tfrac{\mm-\LL}{\nn-\mm}}
	\Quad{given}
		\LL_\Sigma
	=
		\LL.
	\]
\end{cor}

\subsubsection{Justification of Poisson Approximation}
\label{sec:ip:main:just}

It remains to justify why we can replace
$\XX$ and $\YY$
by their Poisson approximations
\[
	\XX_0
\sim
	\Pois(\tfrac \LL\mm)
\Qand
	\YY_0
\sim
	\Pois(1 - \tfrac{\mm-\LL}{\nn-\mm})
\Quad{given}
	\LL_\Sigma
=
	\LL.
\]

\begin{prop}
\label{res:ip:main:just}
	Let $\alpha \in (0, 1)$ and $\kk = \alpha \nn$.
	Let $\LL \sim \HG(\mm, \mm, \nn)$. Let $\XX_0 \sim \Pois(\tfrac \LL\mm)$ and $\YY_0 \sim \Pois(1 - \tfrac{\mm-\LL}{\nn-\mm})$ independently given $\LL$.
	Let $\cc \in \mbr$ and $\tt \cq \tt_\cc$. Write $\beta \cq \ee^{-\cc} / \alpha$.
	Then,
	\[
		\DD1[\tt]
	=
		\exb{ \absb{ \ee^{-\beta} \exb{ (1 + \beta)^{\XX_0 + \YY_0} \mid \LL_\Sigma } - 1 } }
	+	\oh1.
	\]
\end{prop}

\begin{Proof}
\cref{res:ip:main:abst:approx-split} says that
\[
	\DD1[\tt]
=
	\exb{ \absb{ e^{-\beta} \exb{ (1 + \beta)^{\XX+\YY} \mid \LL_\Sigma } - 1 } }.
\]
We look at the difference between this and the replacement of $(\XX, \YY)$ by $(\XX_0, \YY_0)$:
\[
	\Delta
&
\cq
	\absb{ 
		\exb{ \absb{ e^{-\beta} \exb{ (1 + \beta)^{\XX   + \YY  } \mid \LL } - 1 } }
	-	\exb{ \absb{ e^{-\beta} \exb{ (1 + \beta)^{\XX_0 + \YY_0} \mid \LL } - 1 } }
	}
\\&
\le
	\ee^{-\beta}
	\exb{ \absb{ \exb{ (1 + \beta)^{\XX+\YY} - (1 + \beta)^{\XX_0+\YY_0} \mid \LL} } }
\\&
\le
	\ee^{-\beta}
	\exb{ \absb{ (1 + \beta)^{\XX+\YY} - (1 + \beta)^{\XX_0+\YY_0} } }.
\]
Given $\LL$, both $\XX \to^d \XX_0$ and $\YY \to^d \YY_0$.
This is not sufficient a priori since we are looking at the \mgfs---the map $\zz \mapsto (1 + \beta)^\zz$ is unbounded, so we cannot uniformly bound $(1 + \beta)^{\XX + \YY}$.
However, the tails of all of $\bra{ \XX, \YY, \XX_0, \YY_0 }$ decay super-exponentially---as factorials, in fact---whilst $\zz \mapsto (1 + \beta)^\zz$ is only exponential.
Thus, the convergence is sufficiently strong to deduce that \tozero \Delta \asinf \nn.
We give further details, but omit some technical ones at the end.

We use the explicit descriptions for the laws of $\XX$ and $\YY$ given $\LL$ (\cref{res:ip:main:exact:fp-uncond,res:ip:main:exact:fp-cond}) as well as their conditional independence given $\LL$ (\cref{res:ip:main:exact:indep}):
\[
&	\exb{ (1 + \beta)^{\XX + \YY} \mid \LL }
=
	\sumt{\xx,\yy=0}[\infty]
	(1 + \beta)^{\xx+\yy}
	\pr{ \XX = \xx \mid \LL }
	\pr{ \YY = \yy \mid \LL }
\\&\qquad
=
	\sumd{\xx,\yy=0}[\infty]
	(1+\beta)^{\xx+\yy}
\cdot
	\frac1{\xx!}
	\sumd{\ii=0}[\infty]
	\frac{(-1)^\ii}{\ii!}
	\frac{(\LL)_{\xx+\ii}}{(\mm)_{\xx+\ii}}
\cdot
	\frac1{\yy!}
	\sumd{\jj=0}[\infty]
	\frac{(-1)^\jj}{\jj!}
	\frac{(\nn-2\mm+\LL)_{\yy+\jj}}{(\nn-\mm)_{\yy+\jj}}
\\&\qquad
=
	\sumd{\xx,\yy}
	\sumd{\ii,\jj}
	\frac{(-1)^{\ii+\jj}}{\xx! \yy! \ii! \jj!}
	\rbbb{ \frac \LL\mm (1 + \beta) }^\xx
	\rbbb{ \rbbb{ 1 - \frac{\mm-\LL}{\nn-\mm} } (1 + \beta) }^\yy
	\rbbb{ \frac \LL\mm }^\ii
	\rbbb{ 1 - \frac{\mm-\LL}{\nn-\mm} }^\jj
\\&\qquad\hspace*{4em}
\cdot
	\rbbb{ 
		\frac{(\LL)_{\xx+\ii}}{(\mm)_{\ii+\xx}}
	\bigg/
		\rbbb{ \frac \LL\mm }^{\xx+\ii}
	}
\cdot
	\rbbb{ 
		\frac{(\nn-2\mm+\LL)_{\yy+\jj}}{(\nn-\mm)_{\yy+\jj}}
	\bigg/
		\rbbb{ 1 - \frac{\mm-\LL}{\nn-\mm} }^{\yy+\jj}
	},
\]
where the second equality merely rearranges the terms in order to force it to look like
\[
&	\exb{ (1 + \beta)^{\XX_0 + \YY_0} \mid \LL }
=
	\sumt{\xx}
	\tfrac1\xx! \ee^{-\xx} \rbb{ \tfrac \LL\mm \rbr{1 + \beta} }^\xx
	\sumt{\yy}
	\tfrac1\yy! \ee^{-\xx} \rbb{ \rbr{ 1 - \tfrac{\mm-\LL}{\nn-\mm} } \rbr{1+\beta} }^\yy
\\&\qquad
=
	\sumd{\xx,\yy}
	\sumd{\ii,\jj}
	\frac{(-1)^{\ii+\jj}}{\xx! \yy! \ii! \jj!}
	\rbbb{ \frac \LL\mm (1 + \beta) }^\xx
	\rbbb{ \rbbb{ 1 - \frac{\mm-\LL}{\nn-\mm} } (1 + \beta) }^\yy
	\rbbb{ \frac \LL\mm }^\ii
	\rbbb{ 1 - \frac{\mm-\LL}{\nn-\mm} }^\jj.
\]
It remains to control the `error factor' in our expression for
\(
	\exb{ (1 + \beta)^{\XX + \YY} \mid \LL }.
\)
To~this~end,~let
\[
	\delta
&
\cq
	\exb{ (1 + \beta)^{\XX + \YY} \mid \LL }
-	\exb{ (1 + \beta)^{\XX_0 + \YY_0} \mid \LL }
\\&
=
	\sumd{\xx,\yy}
	\sumd{\ii,\jj}
	\frac{(-1)^{\ii+\jj}}{\xx! \yy! \ii! \jj!}
	\rbbb{ \frac \LL\mm (1 + \beta) }^\xx
	\rbbb{ \rbbb{ 1 - \frac{\mm-\LL}{\nn-\mm} } (1 + \beta) }^\yy
	\rbbb{ \frac \LL\mm }^\ii
	\rbbb{ 1 - \frac{\mm-\LL}{\nn-\mm} }^\jj
\\&\qquad\hspace*{3em}
\cdot
	\rbbb{ 
		\rbbb{ 
			\frac{(\LL)_{\xx+\ii}}{(\mm)_{\xx+\ii}}
		\bigg/
			\rbbb{ \frac \LL\mm }^{\xx+\ii}
		}
	\cdot
		\rbbb{ 
			\frac{(\nn-2\mm+\LL)_{\yy+\jj}}{(\nn-\mm)_{\yy+\jj}}
		\bigg/
			\rbbb{ 1 - \frac{\mm-\LL}{\nn-\mm} }^{\yy+\jj}
		}
	-	1
	},
\]
Thus,
\(
	\Delta
\le
	\ee^{-\beta} \ex{ \abs \delta }.
\)
It remains to show that $\ex{ \abs \delta } \to 0$ \asinf \nn to complete the proof.

The ratios of falling factorials are always smaller than their corresponding powers, because $\LL \le \MM$.
Also, $\mm-\LL \le \nn-\mm$, ie $\nn-2\mm+\LL \ge 0$.
Thus,
\begin{gather*}
	\abs \delta
\le
	\sumd{\xx,\yy}
	\sumd{\ii,\jj}
	\frac{(1 + \beta)^{\xx + \yy}}{\xx! \yy! \ii! \jj!}
	\rbb{ 1 - (1 - \aa_{\xx,\ii}) (1 - \bb_{\yy,\jj}) }
\quad
	\text{where}
\\
	\aa_{\xx,\ii}
\cq
	1
-	\frac{(\LL)_{\xx+\ii}}{(\mm)_{\xx+\ii}}
\bigg/
	\rbbb{ \frac \LL\mm }^{\xx+\ii}
\Qand
	\bb_{\yy,\jj}
\cq
	\frac{(\nn-2\mm+\LL)_{\yy+\jj}}{(\nn-\mm)_{\yy+\jj}}
\bigg/
	\rbbb{ 1 - \frac{\mm-\LL}{\nn-\mm} }^{\yy+\jj}.
\end{gather*}
Both $\aa$ and $\bb$ should be viewed as ``approximately $0$''; in particular they are in $[0, 1]$.
Now,
\[
	0
\le
	1 - (1 - \aa)(1 - \bb)
=
	\aa + \bb - \aa \bb
\le
	\aa + \bb
\Qforall
	\aa, \bb \in [0, 1].
\]
This way, we can partially separate the sums over $(\xx,\ii)$ and $(\yy,\jj)$:
\[
	\abs \delta
&
\le
	\sumd{\xx,\ii}
	\dfrac{(1+\beta)^\xx}{\xx! \ii!}
	\aa_{\xx,\ii}
	\sumd{\yy,\jj}
	\dfrac{(1+\beta)^\yy}{\yy! \jj!}
+	\sumd{\yy,\jj}
	\dfrac{(1+\beta)^\yy}{\yy! \jj!}
	\bb_{\yy,\jj}
	\sumd{\xx,\ii}
	\dfrac{(1+\beta)^\xx}{\xx! \ii!}
\\&
\le
	2 \ee^{1+\beta}
	\sumt{\zz,\kk}
	(1+\beta)^\zz (\aa_{\zz,\kk} + \bb_{\zz,\kk}) / (\zz! \kk!).
\]

Algebraic manipulations show that $\abs \delta \to 0$ \asinf \nn when $\min\bra{\LL, \nn-2\mm+\LL} \ge \nn^{1/4}$ for fixed $\beta$, completing the proof.
There are four key facts required for establishing this.
\begin{enumerate}
	\item 
	$\aa_{\zz,\kk}, \bb_{\zz,\kk} \in [0, 1]$ for all $(\zz, \kk)$.
	
	\item 
	The factorial $\zz!$ grows faster than the exponential $(1 + \beta)^\zz$, meaning that asymptotically all the mass of $\sumt{\zz}$ is distributed over $\zz$ with $\zz \asymp 1$.
	
	\item 
	Similarly, asymptotically all the mass of $\sumt{\kk}$ is distributed over $\kk$ with $\kk \asymp 1$.
	
	\item 
	The approximations
		$\aa_{\xx,\ii} \approx 0$ and $\bb_{\yy,\jj} \approx 0$
	hold if
		$(\xx+\ii)^2 \ll \LL$ and $(\yy+\jj)^2 \ll \nn-2\mm+\LL$.
\end{enumerate}
Now, $\LL \sim \HG(\mm, \mm, \nn)$ and $\mm \asymp \nn$, so the condition $\min\bra{\LL, \nn-2\mm+\LL} \ge \nn^{1/4}$ holds with at least exponentially high probability when $\alpha (1 - \alpha) \asymp 1$.
Indeed,
\[
	\ext{ \LL } = \mm^2 / \nn = \alpha^2 \nn
\Qand
	\ext{ \nn - 2\mm + \LL } = \nn - 2\mm + \mm^2 / \nn = (1 - \alpha^2) \nn.
\]
On the other hand, a $\Pois(1)$ \RV is order $\nn$ with super-exponentially low probability; the same holds for $\XX$ and $\YY$.
Thus, $\ex{ \abs \delta } \to 0$ \asinf \nn.
We omit the technical manipulations required to rigorously establish these claims in this last paragraph.
\end{Proof}

\begin{rmkt}
\label{rmk:ip:main:k-small}
We recall that $\beta = \ee^{-\cc} / \alpha$ and $\alpha = \kk/\nn$.
The above proof assumed that $\beta \asymp 1$.
However, all it really needs is for the factorial $1/\rr!$ to beat the exponential $(1 + \beta)^\rr$ in an appropriate range.
If $\beta \asymp 1$, then asymptotically all the mass of the sums is in the first order-$1$ number of terms.
However, if $\beta \gg 1$, then divergently many terms will need to be considered, adding technical difficulties.
This makes some approximations more challenging.

We have not pushed the technical manipulations to see how large $\beta$ can be taken, ie how small $\alpha = \kk/\nn$ can be taken.
We \emph{should} be able to handle \emph{some} $\alpha \ll 1$.
In this case, the limit is
\(
	\pr{ \Pois(\ee^{-\cc}) = 0 }
=
	\exp{-\ee^{-\cc}},
\)
which is the \cdf of the Gumbel distribution.
\end{rmkt}

\paragraph*{Acknowledgements.}

We thank the anonymous referee for their comments, including extra references and suggestions for minor rewriting of statements, as well as a few typographical errors.
The paper is clearer and more readable as a result, for which they have our thanks.

\section*{Bibliography}
\addcontentsline{toc}{section}{Bibliography}

\renewcommand{\bibfont}{\sffamily}
\renewcommand{\bibfont}{\sffamily\small}
\printbibliography[heading=none]

%
%

\end{document}